\documentclass[AMS,Times1COL]{WileyNJDv5} 

\usepackage[printonlyused,withpage]{acronym}
\usepackage{nomencl,framed}
\usepackage{soul}
\usepackage{mathtools}
\usepackage{subcaption}
\articletype{Original Article}%

\received{Date Month Year}
\revised{Date Month Year}
\accepted{Date Month Year}
\journal{Studies in Applied Mathematics}
\volume{00}
\copyyear{2023}
\startpage{1}

 \makenomenclature

\makeatletter
\patchcmd{\thenomenclature}{\section*}{\section}{}{}
\makeatother

\DeclareMathOperator*{\argmin}{arg\,min}
\newcommand{\Mu}{\mathrm{u}}
\newcommand{\Mbu}{\boldsymbol{\mathrm{u}}}
\newcommand{\U}{\Mu^h}
\newcommand{\M}{\mathsf{M}}
 \newcommand{\K}{\mathsf{K}}
\newcommand{\Me}{{e}}
 \newcommand{\Fu}{\mathfrak{u}}
 \newcommand{\Fc}{\mathfrak{c}}
 \newcommand{\D}{\mathsf{D}}

 \newcommand{\Fq}{\mathfrak{q}}
 \newcommand{\Mbf}{\boldsymbol{\mathrm{f}}}
 \newcommand{\HS}[1]{\|#1\|_2}
 \newcommand{\op}[1]{\left\|#1\right\|_\mathrm{op}}
 \newcommand{\BPL}{\mathrm{BPL}}
 \newcommand{\RES}{\mathcal{IR}es}
 \newcommand{\Res}{\mathcal{R}es}
 \newcommand{\up}[1]{$^{\mathrm{#1}}$}

 \newcommand{\Fbu}{\boldsymbol{\mathfrak{u}}}
 \newcommand{\bu}{{{u}}}
 \newcommand\dg\mathfrak
 \newcommand{\B}{\mathcal{B}}
 \newcommand{\A}{\mathcal{A}}
 \newcommand{\mA}{\mathsf{A}}
 
 \newcommand{\MU}{\mathrm{U}}
 \newcommand{\bD}{\boldsymbol{\mathrm{D}}}
 \newcommand{\e}{\mathrm{e}}
 \newcommand{\PS}{*}

 \newcommand{\bU}{{\mathsf{u}}}
 \renewcommand{\d}{\mathrm{d}}

 \newcommand{\T}{\mathsf{T}}
 \newcommand{\slope}{\mathfrak{s}}

\newcommand{\inn}{\textrm{in}}
\newcommand{\rb}{\textrm{b}}

\definecolor{myred}{RGB}{255,0,0}

\title{Stabilized Time Series Expansions for High-Order Finite Element Solutions of Partial Differential Equations}

\author[1]{Ahmad Deeb}

\author[1,2]{Denys Dutykh}

\authormark{DEEB \textsc{et al.}}
\titlemark{STABILIZED TIME SERIES EXPANSIONS FOR HIGH-ORDER FINITE ELEMENT SOLUTIONS OF PARTIAL DIFFERENTIAL EQUATIONS}

\address[1]{\orgdiv{Mathematics Department}, \orgname{Khalifa University of Science and Technology, PO Box 127788}, \orgaddress{\state{Abu Dhabi}, \country{United Arab Emirates}}}

\address[2]{\orgdiv{Causal Dynamics Pty Ltd}, \orgaddress{\state{Perth}, \country{Australia}}}

\corres{Corresponding author Ahmad Deeb. \email{ahmad.deeb@ku.ac.ae}}

\abstract[Abstract]{Over the past decade, \ac{FEM} has served as a foundational numerical framework for approximating the terms of \ac{TSE} as solutions to transient \ac{PDE}. However, the application of high-order \ac{FE} to certain classes of \ac{PDE}s, such as diffusion equations and the \ac{NS} equations, often leads to numerical instabilities. These instabilities limit the number of valid terms in the series, though the efficiency of time series integration even when resummation techniques like the \ac{BPL} integrators are employed. In this study, we introduce a novel variational formulation for computing the terms of a \ac{TSE} associated with a given \ac{PDE} using higher-order \ac{FE}s. Our approach involves the incorporation of artificial diffusion terms on the left-hand side of the equations corresponding to each power in the series, serving as a stabilization technique. We demonstrate that this method can be interpreted as a minimization of an energy functional, wherein the total variations of the unknowns are considered. Furthermore, we establish that the coefficients of the artificial diffusion for each term in the series obey a recurrence relation, which can be determined by minimizing the condition number of the associated linear system. We highlight the link between the proposed technique and the \ac{DMP} of the heat equation. We show, via numerical experiments, how the proposed technique allows having additional valid terms of the series that will be substantial in enlarging the stability domain of the \ac{BPL} integrators.}


\jnlcitation{\cname{%
\author{Deeb A} and
\author{Dutykh D}}.
\ctitle{Stabilized Time Series Expansions for High-Order Finite Element Solutions of Partial Differential Equations.} \cjournal{\it Studies in Applied Mathematics} \cvol{2023;00(00):1--18}.}

\begin{document}
\maketitle



\section{Introduction}\label{sec1}

The numerical integration in time for differential equations has been widely investigated. Numerical discretization schemes were proposed to provide approximation of solutions  \cite{book:hairer,book:hairer2,hairer2002geometric,book:tomas,book:iserles,casas-16}. Technically, classical numerical methods could be classified under two categories:  (1) one step methods where integration of solution for the next time layer is based only on the last approximation, as for instance the celebrated Runge-Kutta schemes \cite{Butcher_1996}, and (2) multi-step methods that use the last $k$ layers to integrate the solution, as for instance the ones based on the \ac{BDF} \cite{book:hairer}. These time numerical methods have been applied to integrate \ac{ODE} and semi-discrete \ac{PDE}s, in addition to numerical methods for the space discretization such as \ac{FVM}, \ac{FDM} and \ac{FEM}.

A new and different type of time numerical integration, that is based on divergent series resummation technique \cite{ramis1993,xups_thomann}, has emerged in the last two decades. The idea was developed for the first time in the PhD thesis of Razafindralandy \cite{dina-thesis}. To start, consider that the solution is written as a time series expansion, and if the series belongs to a categorized class of divergence, namely a Gevrey class \cite{gevrey_1918}, resummation techniques, as Borel-Laplace resummation \cite{borel-1901} could be applied to provide analytical representation of the original series. In his thesis \cite{dina-thesis}, the \ac{BPL} algorithm was developed as a numerical version of the Borel-Laplace resummation method, where the \ac{Pa-ap} was used in the prolongation process.
The scheme was used also to solve classical problems in mechanics and fluid dynamics in the framework of \ac{FDM} for \ac{PDE}s \cite{dina-2012}. Then, it was extended by Deeb \emph{et al.} \cite{deeb-thesis} to solve \ac{PDE}s in the \ac{FEMF}. Other types of problems were also considered by this algorithm as stiff ones \cite{DEEB_2022_bpl} and problems with large time dynamics \cite{ahmad_robust_integrators_2019}.

In this paper we are interested in providing numerical discretizations for transient \ac{PDE}s in the \ac{FEMF} \cite{dougalis-2008} and divergent series resummation technique. After considering that the solution is written in its \ac{TSE} form, terms of the series are space functions and obey a recurrence formula specific to each \ac{PDE}. If the initial condition is an analytical function in space that is explicitly given, then terms are obtained directly by carrying out the recurrence formula to any order.
This is not the case for every \ac{IVP} and even is not practically applicable for many reasons as for instance the following one:
after computing terms up to a fixed order, the solution is approximated by one of the resummation techniques \cite{ramis1993,xups_thomann} where the approximation is valid on a finite interval. To continue evolving, the last approximation is considered to be the initial condition of a modified problem and the above process is then repeated. Within iterations, it is becoming complex and difficult to give explicitly the forms of terms of the series. That is why a numerical approach is followed to deal with the recurrence formula and compute terms in the \ac{FEMF}.

When dealing with \ac{FEM}, it was shown \cite[chapter 5]{deeb-thesis} that for first order \ac{FE}s, terms of the series are computed numerically and their validity is acceptable up to a given order related to the size of mesh. Lumping technique of mass matrix has improved the situation and increased that order.
Nevertheless, first order \ac{FE}s are not sufficient when dealing with different types of problems. For example, solving interaction fluids structure problems in the framework of \ac{LSM} calls for higher order space derivative to compute the effect of bending energy on solid boundaries \cite{deeb:MMA,deeb:casson}. Thus, interpolating of elements should be also done with higher orders too.

Another example is when we deal with the incompressible \ac{NS} equations: in their discrete forms \cite{dougalis-1982}, the $\inf$-$\sup$ condition imposes difference between the orders of \ac{FE}s for velocity and pressure, as for instance the Taylor-Hood element which uses a second order for the velocity and first order for pressure.
In this case, numerical computation of terms of the velocity is not efficient and terms are not valid starting from the second order \ac{FE}.
Applying lumping technique could help improving the computation of terms and hand out a computation of the second term with the second order \ac{FE} but the third term will generally explode. Other techniques such as preconditioning or Multi-Grid methods were tested too but did not give the desirable results.
This is why we are looking for alternative techniques to produce a stable process for computing terms of the \ac{TSE} solution of \ac{PDE}s in a higher-order \ac{FEMF}.

The outline of the paper is as follows. The next Section presents resummation technique we adopt in case of system of \ac{ODE}s that will be useful to semi-discrete forms of PDEs.
Section 3 presents the issue of instabilities we face when solving \ac{PDE}s in a high-order \ac{FEMF}.
Section 4 presents an essay to explain instabilities and how can we limit them.
Section 5 presents the process we apply to overcome these instabilities.
Section 6 presents results of simulation we have done and improvement we have achieved.
We end with main conclusions and perspectives in Section \ref{sec-conc}.

\section{Resummation techniques in case of ODEs}
We consider a system of \ac{ODE}s represented as follows:
\begin{equation}
	\label{prob_i}
	\frac{\d \Mbu}{\d t} = \Mbf(t,\Mbu(t)),\quad t\in\,]t_0,T],
\end{equation}
with a given initial condition $\Mbu_0:= \Mbu(t_0)$. We consider the situation where we have $\Mbu(t) = \big(\Fbu^1(t),\dots,\Fbu^\Fq(t) \big)^{\top}  \in \mathbb{R}^{\Fq}$ is a vector function of $\Fq$ components with $\T$ is the transpose operator.
\begin{equation}
	\Mbu:
	\begin{array}{ccc}
		[t_0,T] &\longrightarrow & \mathbb{R}^\Fq\\
		t & \mapsto & \Mbu(t)
	\end{array},
	\qquad
	\Mbf:
	\begin{array}{ccc}
		[t_0,T]\times \mathbb{R}^\Fq &\longrightarrow & \mathbb{R}^\Fq\\
		\big(t, \Mbu(t)\big) & \mapsto & \Mbf \big(t,\Mbu(t)\big) .
	\end{array}
\end{equation}
Finding the solution to this \ac{IVP} will be based on divergent series resummation. We will consider that the vector solutionc $\Mbu$ and its components $\Fbu^i$ are represented, respectively by $\hat \Fbu$ and $\hat{\Fbu}^i$, via a \ac{TSE} as follows:
\begin{equation}
	\label{serie-rep}
	\hat{\Fbu}^i(t) := \sum\limits_{k \geqslant 0}\Fu^i_k t^k\,, i= 1,\ldots,\Fq\,,\quad \hat{\Mbu}(t) := \sum\limits_{k\geqslant 0}\Mu_k t^k.
\end{equation}
Here, $\Mu_k = \big(\Fu^1_{k},\ldots,\Fu^\Fq_k\big)^{\top}\in \mathbb{R}^\Fq$ are vectors of $\Fq$ components and are governed by a recurrence formula resulting from substituting $\Mbu$ in the system by the series $\hat \Mbu$ and equating coefficients of $t^k$:
\begin{equation}
	\label{recurrence_formula_ode}
	\Mu_{k+1} = F_k(t_0,k,\Mu_0,\ldots,\Mu_k),
\end{equation}
where $F_k$ are polynomial functions. For example, the scalar inhomogeneous linear equation (here $\Fq=1$)  $\displaystyle \frac{\d \Mbu}{\d t} = a\Mbu + \e^{bt}$ has the recurrence formula $\Mu_{k+1} = \displaystyle \frac{1}{k+1}\big(a\Mu_k + \e^{bt_0}\frac{b^k}{k!}\big)$, with $\Mu_0=\Mbu(t_0)$. Having a recurrence formula will give us the ability to compute $\Mu_k$ up to any finite order $m$ and using them to approximate the solution, thanks to \eqref{serie-rep}. The automatic differentiation could be also used to compute the terms of the series \cite{autodiff1,Bucker2005ABo}.

\nomenclature{$t$}{Time variable}
\nomenclature{$t_0$}{Lower bound of time}
\nomenclature{$T$}{Upper bound of time}
\nomenclature{$\Mbu$}{Vector solution of system with $\Fq$ ODEs}
\nomenclature{${\Fbu}^i$}{$i$\up{th} component of vector $\Mbu$}
\nomenclature{$\Mbf$}{Right hand side of the system}
\nomenclature{$\Fq$}{Dimension of the system}
\nomenclature{$\hat{\Mbu},\hat{\Fbu}^i$}{Time series expansion of $\Mbu,\Fbu^i$}
\nomenclature{$\Mu_k$}{$k$\up{th} term in the series $\hat\Mbu$ / Vector element of all values of $u_k$ on nodes $x^i\in \mathcal{T}_h$}
\nomenclature{$\Fu^i_k$}{$k$\up{th} term in the series $\hat\Fbu^i$}

\subsection{Partial sum}
In the case where the series is convergent, the partial sum is a tool that will let us approximate the solution with a given precision $\epsilon$  (see \cite{cochelin_1994}), as follows:
\begin{equation}
	\label{partial_sum}
	\|\Mbu(t) - \sum\limits_{k=0}^{m}\Mu_kt^k\| \leqslant \epsilon \Longleftrightarrow  t \leqslant \sqrt[m-1]{\frac{\epsilon\|\Mu_1\|}{\|\Mu_m\|}},
\end{equation}
where $\|.\|$ is a norm defined in $\mathbb{R}^\Fq$. We denote by $\Delta t^{\PS}$ the greatest $t$ that verifies the left inequality which is also the upper limit of the right one.

\nomenclature{$\epsilon$}{User precision}
\nomenclature{$m$}{Truncated order of the series}
\nomenclature{$\Delta t^*$}{Radius of validity domain of the partial sum}

\subsection{Borel-Padé-Laplace integrator}
However, in the case where the series is divergent, the partial sum becomes inefficient as $\Delta t^{\PS}$ becomes too small to advance simulation in time. This is why we will call for resummation techniques that are able to produce valid approximations for larger $t$. We will employ here the \ac{BPL} integrator \cite{dina-2012,DEEB_2022_bpl} that is based on the Borel-Laplace resummation method \cite{ahmad_bpl_2014}. There are other techniques of resummation, as for instance
the Inverse Factorial Series \cite{ahmad_icnpaa_2016} that is also based on the Borel-Laplace resummation and was compared to the \ac{BPL} algorithm in \cite{ahmad_comp_bpl_sfg_2015}.

\begin{figure*}[!h]
	\centerline{\includegraphics[scale=0.6]{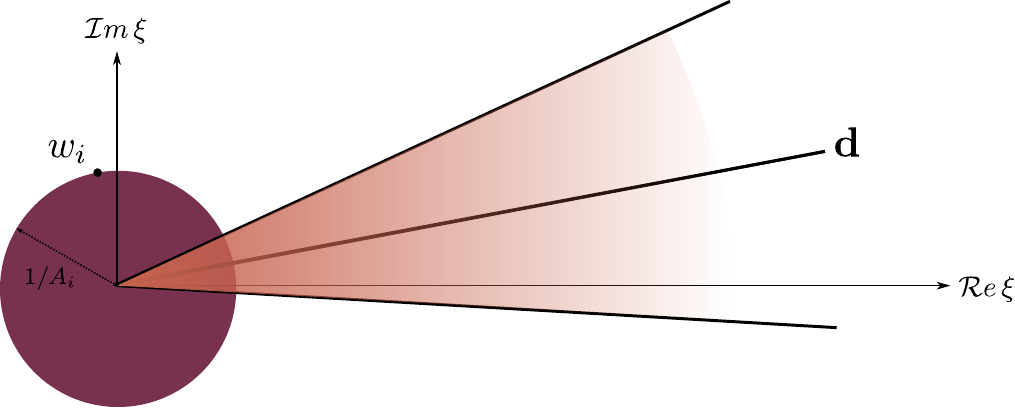}}
	\caption{Illustration, in the Borel space, of convergence domain, singularity and a sector for analytical prolongation.}
	\label{fig_sing}
\end{figure*}

Given a divergent series, we first apply the series version of Borel transformation \cite{balser_1994}, $\hat \B$, on every component $\hat{\Fbu}^i$ in $\hat{\Mbu}$ as follows:
\begin{equation}
	\hat\B(t^k)(\xi) = \frac{\xi^{k-1}}{(k-1)!}.
\end{equation}
Thus, if the series is divergent but not faster than the factorial, \emph{i.e.} there exist positive constants $C_i$ and $A_i$, $i\in\,\{1,\ldots,\Fq\}$, such that for every $k$, we have $|\Fu^i_k|\leqslant C_i A_i^k k!$ ($\Fu^i_k \in \mathbb{R}$), therefore the Borel transformation of the initial series $\hat \Fbu^i$ is a convergent series with $1/A_i$ to be its convergence radius in the complex plane of the Borel space. Despite of having a singularity $w_i$ \cite{polya1929,pauls_frisch_2007} that lies on the circle centred in the origin and radius $1/A_i$, this convergent series could be prolongated analytically along a sector, that does not contain $w_i$, centered by a semi-line $\mathbf{d}$ in the complex plane (see Figure \ref{fig_sing}). For more details, we refer the readers to \cite[Chap. 1]{deeb-thesis}. Numerically, we have $m$ terms of the series, thus the \ac{Pa-ap} could be used to extrapolate numerically the series beyond its convergence domain \cite{pade-brezinski}.

\nomenclature{$\hat\B$}{Series' version of the Borel transformation}
\nomenclature{$\xi$}{Variable in the Borel space}
\nomenclature{$w_i$}{Singularity in the Borel space relative to $\hat \Fbu^i$}
\begin{remark}
	To determine the radius of convergence of $\hat\B \hat{\Mbu}$, we  consider the value $A^{-1}$ to be its convergence radius, where:
	$$A := \sup\limits_{i\in\{1,\ldots,\Fq\} }A_i.$$
\end{remark}

We will use the pointwise \ac{Pa-ap}; \emph{i.e.} for every scalar function $\Fu^i(t), i\in\{1,\ldots,\Fq\}$. There is also another type of \ac{Pa-ap} such as vector Padé that could be considered \cite{pade-app-vecotr-graves}. For a pointwise version and a given series $\hat\Fbu^i$, the \ac{Pa-ap} of order $m-1= r+s$ consists of finding a rational polynomial, denoted by $\dg{P}^i_m(\xi)$:
\begin{equation}
	\dg{P}^i_m(\xi):=\cfrac{\dg a^i_0+ \dg a^i_1\xi+\cdots+\dg a^i_r\xi^r}{\dg b^i_0+\dg b^i_1\xi+\cdots+\dg b^i_s\xi^s}, \quad i\in\,\{1,\dots,\Fq\}
\end{equation}
such that coefficients will verify the following relation:
\begin{equation}
	\hat\Fbu^i(\xi) - \dg P^i_m(\xi) = O(\xi^{m}), \quad i\in\,\{1,\dots,\Fq\}.
\end{equation}
It is possible that the Padé table could not be normal; \emph{i.e.} there is a block in the Padé table where all cases have the same \ac{Pa-ap}. Therefore we use a robust algorithm to compute a \ac{Pa-ap} using the \ac{SVD} technique \cite{pade_svd,ibrayeva-13}. To go back to the physical space, the Laplace transformation $\mathcal{L}$, which is the inverse of the Borel transformation, is applied as follows:
\begin{equation}
	\mathcal{L}\Big(\dg P^i_m(\xi)\Big)(t) = \int_{\mathbf{d}} \dg P^i_m(\xi)e^{-\xi/t} d\xi.
\end{equation}
For simplicity and if there is no singularity that lies on $\mathbf{d}$, we consider the semi-line $\mathbf{d}$ to be the positive real axis $\mathbb{R}^+$, and the Laplace transform could be approximated numerically after a variable change $\mu = -\xi/t$, using the \ac{GL} quadrature:
\begin{equation}
	\int_0^{\infty} \dg P^i_m(\xi)e^{-\xi/t} d\xi = t\int_0^{\infty} \dg P^i_m(\xi t)e^{-\xi} d\xi \approx t\sum\limits_{j=1}^{N_g}\dg{P}^i_m(\xi_j t) \omega_j,
\end{equation}
where $\xi_j$ and $\omega_j$ are respectively $N_g$ points and weights associated to the \ac{GL} quadrature of order $2N_g-1$. The solution is then approximated by the following value at time $t$:
\begin{equation}
	\Fbu^i(t) \approx \Fu^i_0 + \sum\limits_{j=1}^{N_g} \dg{P'}^i_m(\xi_j,t) \,\omega_j, \quad \dg{P'}^i_m(\xi,t) := t\dg{P}^i_m(\xi t).
\end{equation}
To generalize the process and unify the vector notation, we denote by $$\mathrm{P}_m(\xi) := \Big( \dg{P}^1_m(\xi), \dots,\dg{P}^\Fq_m(\xi)\Big)^{\top}$$
the vector of pointwise \ac{Pa-ap}s. Though the vector function $\Mbu(t)$ is approximated by the following:
$$ \Mbu(t) \approx \Mu_0  +  \sum\limits_{j=1}^{N_g}\mathrm{P}'_m(\xi_j,t) \omega_j, \quad \mathrm{P'}_m(\xi,t) = t\mathrm{P}_m(\xi t).$$
To summarize, we present the \ac{BPL} algorithm in Diagram \ref{alg:borel_pade_laplace}.

\captionsetup[table]{name=Diagram,skip=1ex, labelfont=bf}

\begin{table}[!h]
	\caption{\ac{BPL} algorithm.}
	\label{alg:borel_pade_laplace}
	\centering
	{\sf\small
		\begin{tabular}{ccccc}
			$\displaystyle\sum_{k=0}^{m}\Mu_k t^k$& $\dashleftarrow$  &$\displaystyle \Phi^{m}_{t}(\Mu_0) =\Mu_0+\sum_{j=1}^{N_g}\mathrm{P}'^m(\xi_j,t)\,\omega_j$
			\\\\
			$\left.\begin{array}{c}\\\text{Borel}\\\\\end{array}\right\downarrow$&
			&$\left\uparrow \begin{array}{c}\\\text{\ac{GL}} \\\\\end{array}\right.$
			\\\\
			$\displaystyle
			\sum_{k=0}^{m-1}\cfrac{\Mu_{k+1}}{k!}\ \xi^k$
			&$\overrightarrow{\hspace{.5cm}\text{Padé}\hspace{.5cm}}$
			&$\mathrm{P}_m(\xi) = \left( \ldots,\frac{\dg a^i_0+ \dg a^i_1\xi+\cdots+\dg a^i_r\xi^r}{\dg b^i_0+\dg b^i_1\xi+\cdots+\dg b^i_s\xi^s},\ldots \right)^{\top}, {r+s}=m-1$
\end{tabular}}
\end{table}
The dashed left arrow presented in the top of this diagram signifies that if the numerical flow is considered as a function of $t$, then the first $m$ terms of its times series expansion converge to those of the initial series when $N_g\rightarrow \infty$.

\nomenclature{$\dg{P}^i_m$}{\ac{Pa-ap} of order $m$ for series $\hat\Fbu^i$}
\nomenclature{$r$}{Degree of the numerator polynomial in $\dg{P}^i_m$}
\nomenclature{$s$}{Degree of the denominator polynomial in $\dg{P}^i_m$}
\nomenclature{$\mathcal{L}$}{Laplace transformation}
\nomenclature{$N_g$}{Number of \ac{GL} points quadrature}
\nomenclature{$\xi_j,\omega_j$}{$j$\up{th} point and weight of \ac{GL} quadrature}
\nomenclature{$\mathrm{P}_m$}{Vector of pointwise \ac{Pa-ap}}

\captionsetup[table]{name=Table,skip=1ex}

\begin{remark}
It is very important to determine, for given $m,s,r$ and $N_g$ the validity of the approximation for a given $t$; in other words, to estimate the maximum value of $t$ for which the approximation is close to the exact solution $\Mbu(t)$ up to a given precision $\epsilon$, as it was given above for the partial sum. Sure such an error estimation is crucial but it is not available yet, and we will, for moment, skip this part and consider the residual error to estimate the validity of the above approximation; i.e. we consider $\Delta t^{\BPL}$ the maximum value $t$ that verifies the following inequality:
\begin{equation}
	\label{residual-error}
	\left\| \Res(t)(\Mu_0) \right\| = \left\|\frac{\d \Phi^{m}_{t}(\Mu_0)}{\d t} - \Mbf \big(t,\Phi^{m}_{t}(\Mu_0)\big)\right\|\cdot \left\|\frac{\d \Phi^{m}_{t}(\Mu_0)}{\d t} \right\|^{-1} \leqslant \epsilon.
\end{equation}
This issue of determining a priori estimate of $\Delta t^{\BPL}$ will be addressed in future works.
\end{remark}

\subsection{Continuation}
\label{sec_cont}
In this part, we wish to approximate the solution over the segment $[0,T]$. The reader should notice the presence of $\Phi^{m}_{t}$ in the \ac{BPL} diagram: it is the numerical flow resulting from the \ac{BPL} integrator. It is an application that is defined as follows:
\begin{equation}
	\Phi^{m}_{t}
	\left|
	\begin{array}{ccl}
		\mathbb{R}^\Fq& \rightarrow &\mathbb{R}^\Fq\\
		\Mu_0 & \mapsto & \Mu_0+\sum\limits_{j=1}^{N_g}\mathrm{P}'^m(\xi_j,t)\,\omega_j.
	\end{array}
	\right.
\end{equation}
Starting from $t_0 = 0$, we denote by $t_1:=\Delta t^{\BPL}_{1}$ the maximum value of $t$ verifying inequality \eqref{residual-error}. If  $t_1<T$, advancing in time is done by considering a variable change $t = \tau + \Delta t^{\BPL}_{1}$, thus the numerical flow is applied again on a new IVP:
\begin{equation}
	\label{prob_1}
	\frac{\d  \tilde \Mbu_1}{\d\tau} =  \tilde{\Mbf}(\tau,{ \tilde \Mbu_1}(\tau)), \quad \tau \in \,]0,T-t_1],\,  \quad { \tilde \Mbu_1}(0) =\Phi^{m}_{t_1}(\Mbu_0):= \Mu_{0,1}.
\end{equation}
Here $\tilde{\Mbf}\Big(\tau,{ \tilde \Mbu_1}(\tau)\Big) =  \Mbf\Big(\tau+t_1,{ \tilde \Mbu_1}(\tau)\Big) $, therefore we compute again $m$ terms of the new series
\begin{equation}
	\hat{  \tilde\Mbu}_1(\tau) = \sum\limits_{k\geqslant 0} \Mu_{k,1} \tau^k,
\end{equation}
by substituting it in Equation \eqref{prob_1} and identifying coefficients if front of $\tau^k$ to get the recurrence formula for $\Mu_{k+1,1}$:
\begin{equation}
	\Mu_{k+1,1} = F_k(t_1,k,\Mu_{0,1},\ldots,\Mu_{k,1}).
\end{equation}
The numerical flow $\Phi^m_{\tau}$ is then applied for $\tau = \Delta t^{\BPL}_{2}$ to approximate the solution to problem \eqref{prob_i} on $t_2 = t_1 +  \Delta t^{\BPL}_{2}$:
\begin{figure}[!h]
	\centerline{
		\includegraphics[width=0.65\textwidth]{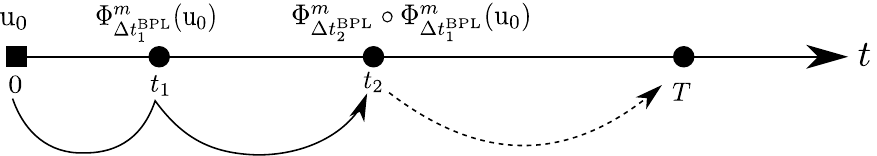}
	}
	\caption{Illustration of the continuation process.}
	\label{}
\end{figure}
\begin{remark}
	If System \eqref{prob_i} is autonomous, then the recurrence formula relative to $\Mu_{k,n}$, terms of the series of the modified problem, is the same formula \eqref{recurrence_formula_ode} of that of $\Mu_k$.
\end{remark}

We denote by $\Mbu_0:=\Mu_{0,0}\equiv \Mu_0$. At the iteration $n$, we consider having $\Mbu_n:=\Phi^m_{\Delta t^{\BPL}_n}(\Mbu_{n-1})$ approximation $\Mbu(t_n = t_{n-1} + \Delta t^{\BPL}_n)$. To find the maximum value $\Delta t_{n+1}^{\BPL}$, we start with $\Delta t^{\BPL}_{n+1}=\Delta t^*$ as in \eqref{partial_sum}, then we compute $\Phi^m_{\Delta t_{n+1}^{\BPL}}(\Mu_{0,n})$ and its residual error. If the latter is less then the pre-defined user tolerance, we increase $\Delta t_{n+1}^{\BPL}$ by $10 \%$ ($\times 1.1$). Algorithm \ref{alg:1} presents the continuation process with these steps.

\begin{algorithm}
	\caption{Continuation with \ac{BPL}}\label{alg:1}
	\begin{algorithmic}
		\Require $t_0$,  $\Mu_0$, $m$, $\epsilon$
		\State $t_n\gets t_0$
		\State $\Mbu_n\gets \Mu_0 $
		\While {$t_n <T $}
		\State $\Mu_{0,n}\gets \Mbu_n $
		\For {$k\gets 0$ to $m-1$}
		\State $\Mu_{k+1,n}\gets F_k(t_n,k,\Mu_{0,n},\ldots,\Mu_{k,n}) $
		\EndFor
		\State $\Delta t_{n+1}^{\BPL} \gets \displaystyle\sqrt[m-1]{\frac{\epsilon\|\Mu_{1,n}\|}{\|\Mu_{m,n}\|}}$
		\ \State $ \Mbu_{n} \gets \Phi^m_{\Delta t^{\BPL}_{n+1}}(\Mu_{0,n})$
		\While{$\left\| \Res( \Delta t^{\BPL}_{n+1})(\Mu_{0,n}) \right\|\leqslant \epsilon$}
		\State
		\State $\Delta t_{n+1}^{\BPL}\gets\Delta t_{n+1}^{\BPL} \times 1.1$
		\State $ \Mbu_{n} \gets \Phi^m_{\Delta t^{\BPL}_{n+1}}(\Mu_{0,n})$
		\EndWhile
		\State $t_{n}\gets t_n + \Delta t^{\BPL}_{n+1}$
		\EndWhile
	\end{algorithmic}
\end{algorithm}

By following the continuation process, we will be able to reach the simulation at the end of the required interval after a finite number of steps. In the next section, we present the problem we have when using time series expansion for some classes of \ac{PDE}s in the \ac{FEMF}.

\nomenclature{$\Phi^m$}{Numerical flow of the \ac{BPL} algorithm}
\nomenclature{$t_n$}{$n$\up{th} point of time}
\nomenclature{$\Mbu_n$}{$n$\up{th} approximation at $t=t_n$}

\section{Problem in FEM discretization of parabolic \ac{PDE}s}

We consider a scalar dependent variable $\bu$
\begin{equation*}
	\bu :
	\begin{array}{ccc}
		[t_0,T] \times \Omega & \rightarrow & \mathbb{R}\\
		(t,x) & \mapsto & \bu(t,x)
	\end{array}
\end{equation*}
that is governed by the following partial differential equation:
\begin{equation*}
	\begin{array}{rcl}
		\partial_t \bu &=& \A(x,\bu,\nabla \bu,...), \quad (t,x) \in (t_0,T]\times \Omega,\\
		\bu(t_0,x) &=& u_0(x), \qquad\qquad x \in \Omega ,
	\end{array}
\end{equation*}
where $t$ and $x$ are the time and space variables respectively, $\Omega$ is the space domain, $\partial_t$ is the partial derivative with respect to time, $\A$ is a differential operator and $u_0(x)$ is the initial condition.
Considering $t_0 = 0$, the solution could be written as times series expansion:
\begin{equation}
	\label{TSE_PDE}
	\hat \bu(t,x) := \sum\limits_{k=0}^{\infty} u_k(x)\, t^k,
\end{equation}
where $u_k(x)$ are space dependent functions obeying a recurrence formula. After injecting Formula \eqref{TSE_PDE} in the last \ac{PDE} and equating coefficients of $t^k$, it will provide us a series of equalities that compute $u_k(x)$ as follows:
\begin{equation}
	\label{recurrence_formula_pde}
	u_{k+1}(x) = \frac{1}{k+1} \A_k\big(t_0,x,u_0(x),u_1(x),\ldots,u_k(x)\big),
\end{equation}
where $\A_k$ contains differential operations.

Numerically, we will use the \ac{FEMF} to approximate differential operators, for what it can deal with complex geometries on one hand and the flexibility of increasing the order of approximations in space on the other hand \cite{dougalis-1996}. We will present the problem we face in a case study of the heat equation with constant diffusion coefficient, for what it is the simplest transient parabolic \ac{PDE} we know that already highlights the numerical instability issue.

\subsection{FEMF}
Consider $\Omega \subset \mathbb{R}^d (d=1,2,3)$ an open and compact domain. After choosing a finite number of nodes (vertices) $x^i\in\Omega, i=1,\ldots,N_x$, we define a triangulation $\mathcal{T}_h$ of the domain $\Omega$. It consists of having $\mathcal{K}_j$ open space elements that verify the following conditions:
\begin{equation}
	\bigcup\limits_{j} \mathcal{K}_j = \Omega, \quad \mathcal{K}_j \cap\mathcal{K}_{j^{\prime}} = \varnothing, \quad \overline{\mathcal{K}_j} \cap\overline{\mathcal{K}_{j^{\prime}}} \subset \mathbb{R}^{d^\prime},\,{d^\prime<d}.
\end{equation}
In practice, we choose elements $\mathcal{K}_j$ to be conformal simplices in $\mathbb{R}^d$; lines in $\mathbb{R}$, triangles in $\mathbb{R}^2$ and tetrahedra in $\mathbb{R}^3$. The first relation shows that the domain $\Omega$ is covered by the unions of all elements $\mathcal{K}_j$, while the second one shows that all elements are disjoint. The third relation shows that elements are connected by points (if $d=1$), or lines and points (if $d=2$) or points, lines and triangles (if $d=3$). Thus, every element has a finite number of vertices $x^i$, that are numbered locally, and form the connections. Though we have the connectivity matrix that relates the global numbering of vertices by their local numbering related to each element $\mathcal{K}_j$. To end this notation part, $h$ is the maximum diameter of circumscribed circles of all elements $\mathcal{K}_j$.

After having a defined triangulation $\mathcal{T}_h$, we will construct the space of functions that will be used to approximate solutions $u_k(x)$, $k\in\,\{1,\ldots,m\}$. For every degree $p$, we define the following space:
\begin{equation}
	\mathbb{V}_{h,p} = \left\lbrace  v\in \mathrm{C}^0(\Omega)\,\mid \,  v\rvert_{\mathcal{K}_j} \in \mathbb{P}_p[x], \forall \mathcal{K}_j \in \mathcal{T}_h \right\rbrace,
\end{equation}
where $\mathbb{P}_p[x]$ denotes the piecewise continuous Lagrange polynomials of order $p$ over $\Omega$. Having, on every space element $\mathcal{K}_j$ with $q$ nodes, nodal values $\Fu^i_{k} := u_k(x^i)$ and basis functions $\phi^i$ ($\phi^i(x^l)  = \delta_{il}$), we construct a continuous approximation of $u_k(x)$ over $\mathcal{K}_j$:
\begin{equation}
	\label{approx_uk_Kj}
	u_k(x)|_{\mathcal{K}_j} \approx   \sum\limits_{x^i \in \mathcal{K}_j}\Fu^i_{k}\phi^i(x), \quad \forall x\in \mathcal{K}_j.
\end{equation}
Assembling approximations over all elements $\mathcal{K}_j$ will lead us to have an approximation of $u_k(x)$, denoted by $u^h_k(x)$:
\begin{equation}
	\label{approx_uk_omega}
	u_k(x) \approx u^h_k(x) := \sum_{i=1}^{N_x} \Fu^i_k \phi^i(x).
\end{equation}
We denote also by $\bu^h(t,x)$ the approximation of $\bu(t,x)$ in the constructed framework:
\begin{equation*}
	\bu(t,x) \approx \bu^h(t,x) := \sum_{i=1}^{N_x} \Fbu^i(t) \phi^i(x),
\end{equation*}
where $\Fbu^i$ is represented by its series expansion that is given below by:
\begin{equation*}
	\hat\Fbu^i(t) \equiv \sum\limits_{k=0}^{\infty} \Fu^i_k t^k.
\end{equation*}
We denote by $\Mbu(t) = \big(\Fbu^1(t),\ldots,\Fbu^i(t),\ldots, \Fbu^{N_x}(t) \big)^{\top}$ the unknown vector function that represents approximation of time evolution of the dependent variable $u$ at all nodes $x^i$, and by  $\Mu_k = (\Fu_k^1, \ldots, \Fu_k^{N_x})$ vectors of unknowns we are looking to approximate in the \ac{FEMF}. To this end, The FEniCs project \cite{LoggEtal2012,ScroggsEtal2022} is used under Python programming language as a representative generic \ac{FEMF}.

\nomenclature{$x$}{Space variable}
\nomenclature{$u$}{Scalar variable depending of $t$ and $x$}
\nomenclature{$u_0(x)$}{Initial condition}
\nomenclature{$u_k(x)$}{$k$\up{th} term of the series $\hat u$}
\nomenclature{$\Omega$}{Space domain}
\nomenclature{$\A$}{Right hand side of \ac{PDE}}
\nomenclature{$\A_k$}{$k$\up{th} term in the Taylor expansion of $\A$}
\nomenclature{$\mA_k$}{Discrete form of $\A_k$}
\nomenclature{$\mathcal{T}_h$}{Triangulation of $\Omega$}
\nomenclature{$\mathcal{K}_j$}{Simplex element of triangulation $\mathcal{T}_j$}
\nomenclature{$h$}{Maximum of all diameters of circumscribed circles to elements $\mathcal{K}_j$}
\nomenclature{$p$}{Polynomial degree}
\nomenclature{$\mathbb{V}_{h,p}$}{Finite Element space}
\nomenclature{$(.)^\T$}{Transpose of a vector}
\nomenclature{$\phi^l$}{Basis element of $\mathbb{V}_{h,p}$}
\nomenclature{$\U_k$}{Approximation of $\Mu_k$ in $\mathbb{V}_{h,p}$ having exact $\Mu_0,\ldots,\Mu_{k-1}$}

\subsection{Case study: heat equation}
\label{sec_case-study}
Take for instance the following example of the one-dimensional heat equation:
\begin{equation*}
	\partial_t \bu = \nu \Delta \bu, \quad (t,x) \in ]0,T[\times \Omega
\end{equation*}
with the initial condition $u_0(x)$, $\Omega\subset \mathbb{R}$, $\nu$ is the diffusion coefficient and $\Delta := \partial_x^2$ is the one dimensional Laplace operator. Inserting the formula \eqref{TSE_PDE} will lead us to a cascade of equalities:
\begin{equation}
	\label{strong-form}
	u_{k+1}(x) = \frac{\nu}{k+1}\Delta u_k(x) = \frac{\nu^{2}}{(k+1)k} \Delta^{2} u_{k-1}(x)= \ldots =\frac{\nu^{k+1}}{(k+1)!} \Delta^{k
		+1} u_0(x).
\end{equation}
In this example, we have an explicit formula for terms $u_k(x)$ $\forall$ $k\in \mathbb{N}$ in function of the initial condition $u_0(x)$, where not only algebraic operations are involved, but also differential operators. The solution was studied by Balser \cite{balser-99} using divergent series, where it was shown that the class of divergence depends on the initial condition $u_0(x)$.

Numerically, we consider a \ac{FEMF} for approximations of $u_k(x)$, $k\in\{1,\ldots, m\}$.
The goal is to find, for a given triangulation $\mathcal{T}_h$ and a space function $\mathbb{V}_{h,p}$ reliable approximations $u^h_{k} \in \mathbb{V}_{h,p}$ for $u_k(x)$.
To do that, we first present the variational formulation of the recurrence formula relative to the problem above. After applying the Green formula and having homogeneous Dirichlet boundary conditions, it leads us to the following weak formulation:
\begin{equation*}
	\int_{\Omega} u_{k+1}v\,dx = -\frac{\nu}{k+1}\int_{\Omega} \nabla u_k \nabla v\,dx, \quad \forall k>0,
\end{equation*}
which should be verified for every test function $v$ belonging to a suitable space function $V$. Considering $V = \mathbb{V}_{h,p}$, and verifying the last relation for every basis function $\phi^i \in  \mathbb{V}_{h,p}$, we have:
\begin{equation*}
	\displaystyle\int_{\mathcal{K}_j} u^h_{k+1}\phi^i\,dx = -\frac{\nu}{k+1}\int_{\mathcal{K}_j} \nabla u^h_{k} \nabla \phi^i\,dx, \quad \forall \, \phi^i, \,\forall\,\mathcal{K}_j \in \mathcal{T}_h.
\end{equation*}
After replacing $u_k^h$ by its formula \eqref{approx_uk_omega} and assembling all elements in the triangulation $\mathcal{T}_h$, we reach the following linear system with respect to $ u^h_{k+1}$:
\begin{equation}
	\label{lin_sys_heat}
	\M^h \U_{k+1} = -\frac{\nu}{k+1}\K^h \U_{k},
\end{equation}
where $\M^h$ and $\K^h$ are mass and stiffness two-rank tensors associated to a given triangulation and space function $\mathbb{V}_{h,p}$:
\begin{align}
	\label{mass_stiff_matrix}
	\M_{i,j}^h &=  \int_{\Omega} \phi^i\phi^j,
	&
	\K_{i,j}^h &= \int_{\Omega} \nabla\phi^i \nabla \phi^j
\end{align}
and $\U_{k} :=\big( \Fu_k^{h,1},\ldots,\Fu_k^{h,N_x}\big)^{\top}$ is a vector sorting all approximations $\Fu_k^{h,i}\approx\Fu_k^i$ for all points $x^i$ in $\mathcal{T}_h$. Solving the above linear system will produce a vector $\Mu_k^{h} $ approximating $\Mu_k$, where we denote by $\Me_k^{h}$ the following absolute error:
\begin{equation*}
	\label{error_vev_uk}
	\Me^{h}_k = \|\Mu_k - \Mu^h_k\|.
\end{equation*}

\begin{figure}[!h]
	\centerline{
		\includegraphics[scale=0.45]{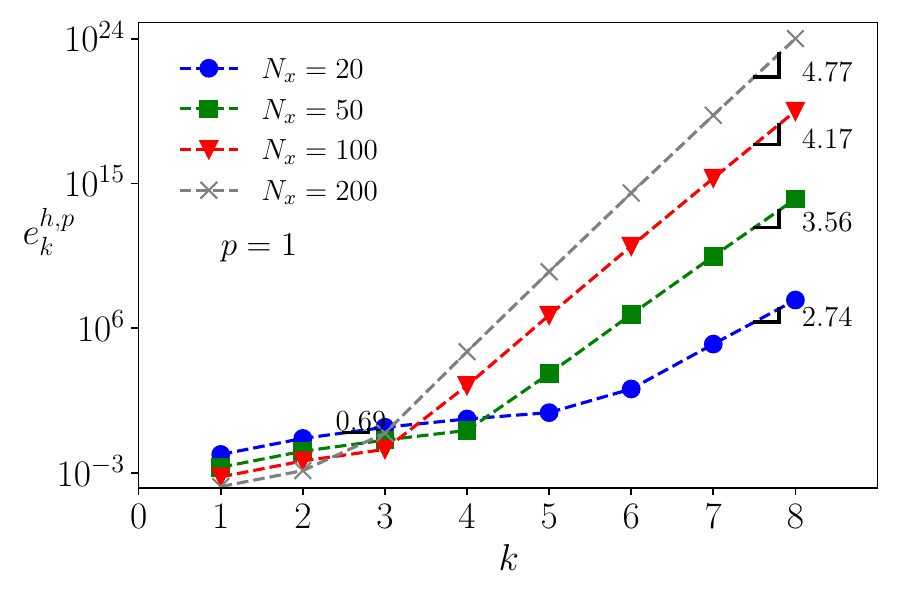}
		\includegraphics[scale=0.45]{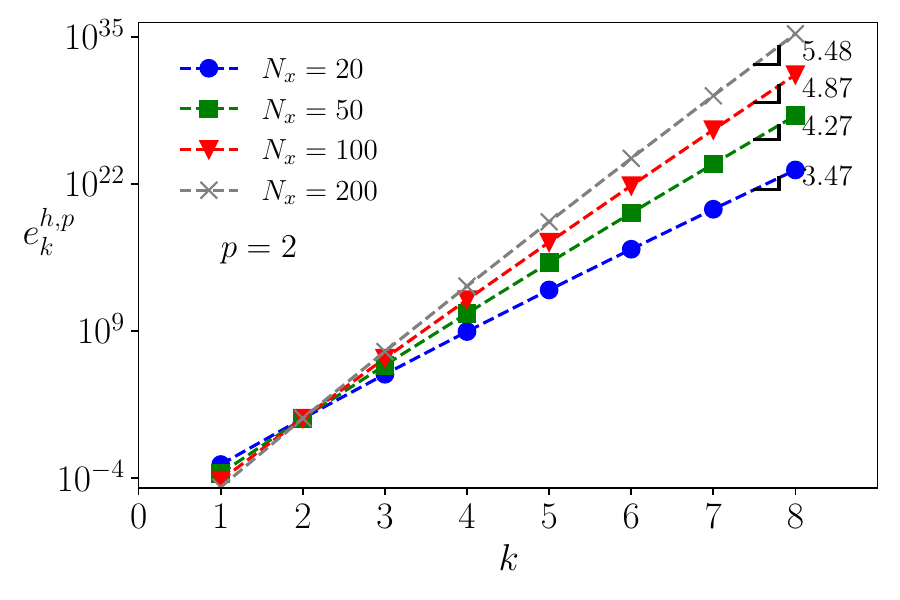}}
	\centerline{
		\includegraphics[scale=0.45]{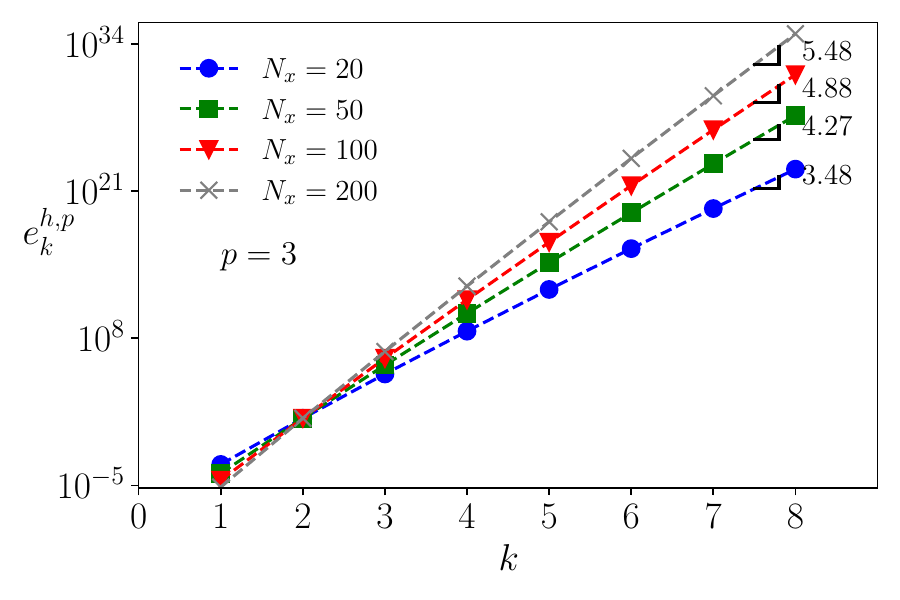}
		\includegraphics[scale=0.45]{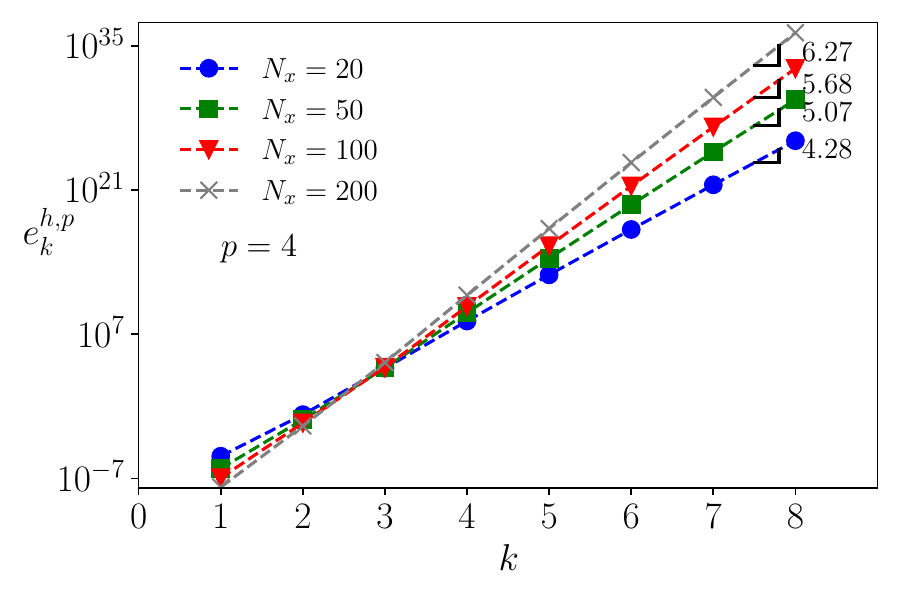}
	}
	\caption{Absolute errors and their slopes $\slope({h,p})$, in the logarithmic scale, for different order $p$ and different mesh sizes for the specific case}
	\label{fig1}
\end{figure}

\paragraph*{Illustration}Consider the case when $\Omega = ]0,1[ $, $u_0(x) = \sin(\pi x)$, thus $u_k(x) = \displaystyle \frac{(-\nu\pi^2)^k}{k!}\sin(\pi x)$ for homogeneous Dirichlet boundary conditions.
By using a first, second, third and fourth order of finite elements, we solve the system \eqref{lin_sys_heat} for a finite number of $k$, and present in Figure \ref{fig1} the absolute errors $\Me^h_k$ for different mesh sizes. The slopes of lines, denoted as $\slope(h,p)$ and defined by the logarithmic errors in the base ten are presented in Figure \ref{fig1} too.
We can see that, when $p=1$, errors $\Me^h_k$ increase within $k$, until exploding at some order that decreases when $h$ decreases. In this case, we can maximize the benefits of resummation technique by lowering the precision of the space $\mathbb{V}_{h,p}$, which is contrary to the high precision in time. Nevertheless, for $p=2,3,4$, the first term is the only one that is computed with acceptable precision and its error decreases when $h$ decreases or $p$ increases. However, errors in terms $\Mu_k, k\geqslant 2$ increases exponentially with $k$, and could be fitted by the following formula:
\begin{equation}
	\label{error_slope}
	\Me_k^{h,p} \sim
	\mathrm{C}(h,p) \times 10^{\,\slope({h,p})\times k},  \quad p\in \{2,3,4\},
\end{equation}
where slopes $\slope(h,p)$ are reported in Table \ref{tab1} and $\mathrm{C}(h,p)$ are fitted constants depending in the \ac{FEMF} and will be determined later on.
\begin{table}[!h]
	\caption{The slopes $\slope({h,p})$ as obtained by Formula \eqref{error_slope} for different \ac{FEMF}.}
	\label{tab1}
	\begin{center}
		\begin{tabular*}{\textwidth}{@{\extracolsep\fill}cccccc@{}}
			\toprule
			 & \multicolumn{2}{@{}|c|}{$p=1$} & $p=2$ & $p=3$ & $p=4$ \\
			\midrule
			$h=1/20$ & 0.69 ($k\leqslant 5$) &2.74($k> 5$) & 3.47 & 3.48 & 4.28 \\
			$h=1/50$  & 0.69 ($k\leqslant 4$)& 3.56($k>4$) & 4.27 & 4.27 & 5.07 \\
			$h=1/100$ & 0.69 ($k\leqslant 3$)& 4.17($k>3$) & 4.87 & 4.88 & 5.68\\
			$h=1/200$ & 0.69 ($k\leqslant 2$)& 4.77($k>2$) & 5.48 & 5.48 & 6.27\\
			\bottomrule
		\end{tabular*}
	\end{center}
\end{table}

In the next Section, we present an attempt to understand the error propagation and its amplification along computing terms $\Mu^h_k$, though we can propose a remedy or alternative process to reduce error amplification.

\nomenclature{$\nu$}{Diffusion coefficient}
\nomenclature{$\M^h$}{Mass Matrix relative to $\mathbb{V}_{h,p}$}
\nomenclature{$\K^h$}{Stiffness matrix relative to $\mathbb{V}_{h,p}$}
\nomenclature{$\Me^h_k$}{Error between the exact vector $\Mu_k$ and the approximation $\bU_k$}
\nomenclature{$N_x$}{Number of nodes $x^i\in \mathcal{T}_h$}

\subsection{Improving the stability with mass lumping}

\begin{figure}[!ht]
  \centerline{
  \includegraphics[scale=0.45]{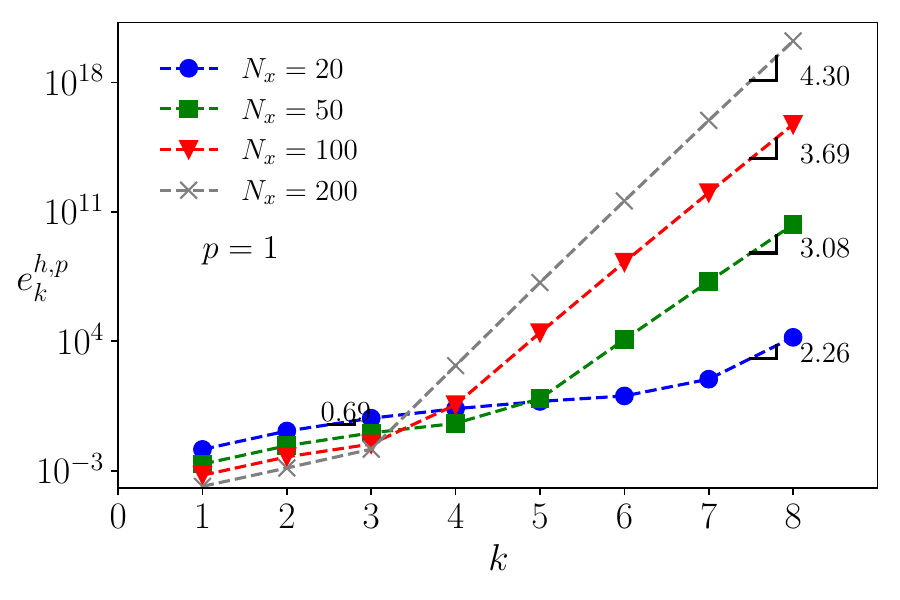}\includegraphics[scale=0.45]{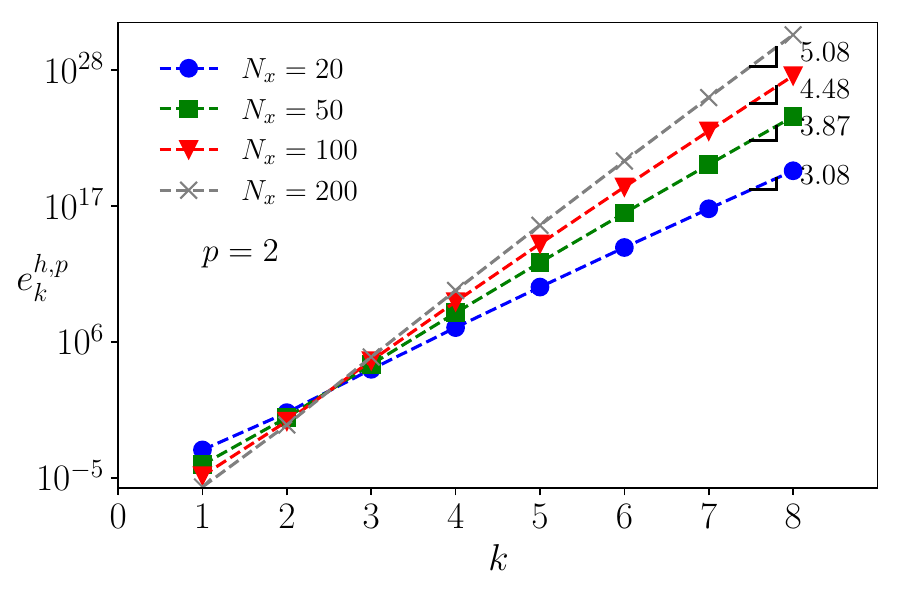}}
\caption{The errors in the terms of the series with lumping technique for the first and second order FE, before stabilization technique.\label{fig_error_lump1}}
\end{figure}

The technique of mass lumping improves the stability in the errors of the terms of the series as was shown in \cite[chapter 5, page 94]{deeb-thesis} in the case of the first order \ac{FE}. Here, the technique of lumping consists of adding the elemetns in every raw of the mass matrix and put the result in the diagonal part of each raw. The left panel in Figure \ref{fig_error_lump1} presents the errors in the terms of the series before the stabilization but with the application of lumping technique for the first order \ac{FE}.
The results show an improvement in the stability as the lumping allows having at least one additional valid term of the serie. In addition, the errors of the terms in the series decrease as shown in the value of the slope for the given mesh sizes. However, the technique of mass lumping for higher order is more complicated. In the following test, we use the Gau\ss-Lobatto integration quadrature to assemble the Mass matrix in the framework of the second order finite element. This quadratur provides also a diagonal mass matrix. The tests we have conducted show only reduction in the slope of the errors as presented in the right panel of Figure \ref{fig_error_lump1}.

\section{Understanding the error amplification}
First, we study the convergence error when computing $\Mu^h_k$, having exact solution of $\Mu_0,\ldots,\Mu_{k-1}$, regarding the size of mesh $h$ and for different degrees $p$. We aim in this part to see the order of the error in $\Mu^h_k$ that will propagate to $\Mu^h_{k+1}$.
To better understand the error propagation, we check the impact of conditioning of matrix $\M^h$, for different degrees $p$, on solving a linear system when the right hand side is perturbed. We end this section by presenting a thorough analysis of how the error, due to discrete form of the equation of $\Mu_k$, will propagate to the approximation $\Mu^h_{k+1}$ and give some conditions to have a stable computation of the series terms.

\subsection{Convergence error in $\Mu^h_1$}
\begin{center}
\begin{figure}[h]
	\centerline{
		\rotatebox[origin=b]{-90}
		{\includegraphics[height=0.45\textheight]{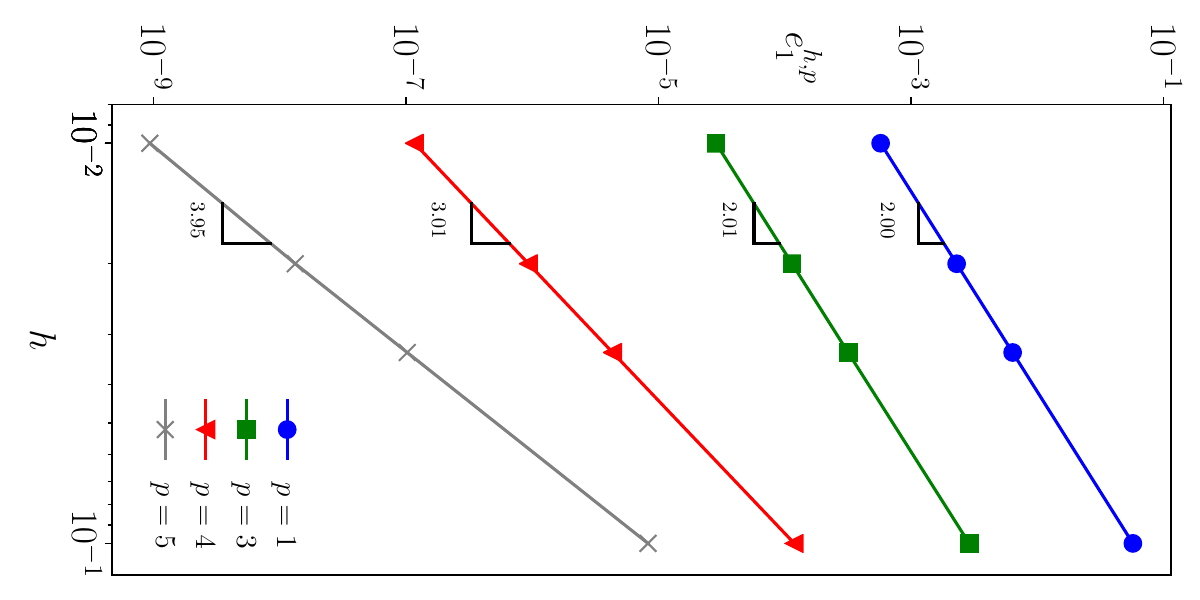}}
	}
	\caption{Plot of the convergence error $\Me^{h,p}_1$ with respect to $h$ and for different degrees $p$.}
	\label{fig_conv}
\end{figure}
\end{center}

We solve the linear system relative to $\Mu^h_1$ for different degrees $p$ and different values $h$, having the initial condition. We denote by $\delta^h \Mu_k$ the following discrepancy vector:
\begin{equation*}
	\delta^h \Mu_k = \U_k - \Mu_k\quad  \text{with}\quad \Me^h_k= \|\delta^h \Mu_k\|.
\end{equation*}
It represents the error in $\Mu_k$ related only to the discretization of the equation, having exact solutions of all previous terms $\Mu_{k-1},\ldots, \Mu_0$. Figure \ref{fig_conv} presents this error for $k=1$. We have inserted the slope of the logarithmic error convergence as it was fitted by the below empirical formula formula:
\begin{equation}
	\label{empirical-form-e_1h}
	\Me^{h,p}_1 \sim
	\left\lbrace
	\begin{array}{rc}
		8^{3/2}\left(\frac{1}{8}\right)^{p-1}\left( \frac{h}{2}\right)^2, & p=1,2,\\
		\left(\frac{h}{2}\right)^{p-1},& p=3,4,5  .
	\end{array}
	\right.
\end{equation}
This formula is valid for any $\Me^h_k$ when having exact solutions for all previous terms $\Mu_{k-1}, \ldots, \Mu_0$. Next, we check the hypothesis of the condition number of mass matrix $\M^h$ and see if there is a correlation between this number and the exploding error in $\Mu_k$ approximation.

\subsection{Condition number of mass matrix}
Formula \eqref{empirical-form-e_1h} shows us that the precision of the approximation $\Mu^h_1$ increases when $p$ increases and $h$ decreases.
Using Figure \ref{fig1}, we can consider that, for $p\geqslant 2$, $\Me^{h,p}_1$ in Equation \eqref{empirical-form-e_1h} is a fitted one for constants $\mathrm{C}(h,p)$ in Equation \eqref{error_slope}:
\begin{equation*}
	\Me^{h,p}_k = \Me^{h,p}_1\times 10^{\slope(h,p)\times (k-1)}, \quad p\geqslant 2,
\end{equation*}
where $\slope(h,p)$ can be retrieved form Table \ref{tab1}. Using this Table, we can remark that the smaller the mesh size $h$ is the bigger the error becomes.
On the other side, we remark that the condition number of $\M^h$, denoted by $\kappa(\M^h)$ and computed for several values of $h$ and $p$, increases when $h$ decreases or when $p$ increases, as shown in Table \ref{tab2}.
\begin{table}[!h]
	\caption{Condition number of matrix $\M^{h}$ for different mesh sizes $h$ and degrees $p$.}
	\label{tab2}
	\begin{center}
		\begin{tabular*}{\textwidth}{@{\extracolsep\fill}lllll@{}}
			\toprule
			& $h=1/10$ & $h=1/30$ & $h=1/50$ & $h=1/100$ \\ \midrule
			$p=1$ &28.60&   89.51&  149.70&  299.85\\
			$p=2$ &63.63 & 195.22 &  325.97&  652.44\\
			$p=3$ &106.87 &  325.48&  543.14 & 1086.85\\
			$p=4$ &158.52 &  480.73 &  801.92 & 1604.43 \\
			\bottomrule
		\end{tabular*}
	\end{center}
\end{table}
It is possible to fit $\kappa(\M^h)$ with the following relation:
\begin{eqnarray}
	\label{Cond_M_formula}
	\kappa(\M^h)\sim \frac{1}{h} \left(\Fc(p) \right)^p, 
\end{eqnarray}
where $\Fc(p)$ is a decreasing function tends to 1 and having $\Fc(1) \approx 3, \Fc(4)\approx 2$.

The condition number $\kappa(\M^h)$ has a direct impact on solving linear system when the right-hand side is perturbed. One has the following linear system:
\begin{equation}
	\label{lin_sys}
	\M \Mu = a
\end{equation}
and wants to see how the solution depends on a perturbed $a$ by $\delta a$, \emph{i.e.} estimate $\delta \Mu$ for a given $\delta a$ as follows:
\begin{equation}
	\M(\Mu+\delta \Mu) = a + \delta a.
\end{equation}
We can show that
\begin{equation}
	\label{eq-cond-1}
	\frac{\|\delta \Mu\|}{\|\Mu\|} \leqslant \kappa(\M)\frac{\|\delta a\|}{\|a\|},
\end{equation}
where $\|.\|$ is the norm of given vectors and $\kappa(\M)$ is the condition number of matrix $\M$ given by:
\begin{equation}
	\kappa(\M) = \op{\M} \cdot \op{ \M^{-1}} \leq   \HS{\M}\cdot\HS{ \M^{-1}}
\end{equation}
with $\op{\M} = \max\limits_{\|\Mu\|=1}\|\M\Mu\|$ is the operator norm and $\HS{\M} = \left[\sum\limits_{i,j}|\M_{i,j}|^2\right]^{1/2} $ is the Frobenius norm of the matrix. From Formula \eqref{eq-cond-1}, we can conclude that the quantity $\|\delta \Mu\| / \|\Mu\|$ resulting from the relative change in the right hand side is limited by the latter multiplied by the condition number, \emph{i.e.} the condition number is the relative error magnification factor.

\nomenclature{$\delta^h \Mu_k$}{Error related to the discrete form of the equation}
\nomenclature{$\slope$}{Slope of line in the logarithmic scale}
\nomenclature{$\kappa(\M)$}{Condition number of matrix $\M$}

\subsection{Error propagation along $\Mu_k$}
It was shown that there is a correlation between the condition number of $\M^h$ and the explosion in error along $\U_k$ when $h$ decreases or $p$ increases: since the approximation of first term $u_1$ presents an error, this error is amplified by the condition number within the approximation of the second term, and so on (\emph{c.f.} Table \ref{tab2}). To see in details how this propagation spreads to the next term, we will consider a recurrence formula relative to a given \ac{PDE}:
\begin{equation}
	u_{k}(x) = \frac{1}{k}\A_{k-1}(t_0,x,u_0,\ldots,u_{k-1}).
\end{equation}
We denote by $\U_{k}$ the approximation of vector $\Mu_k$ representing all values $\Fu_k^i$ of $u_{k}(x)$, obtained in the \ac{FEMF}, at node $x^i$ if all $\Mu_0,\ldots, \Mu_{k-1}$ are known exactly. In this case, the discrete form of the above equation is written below:
\begin{equation}
	\label{sys_exct_un}
	\M^h \U_{k} = \frac{1}{k}\mA_{k-1}(\Mu_0,\ldots,\Mu_{k-1}) = a_{k-1},
\end{equation}
where $\mA_{k-1}$ is the discrete form of $\A_{k-1}$ in the \ac{FEMF}:
\begin{equation}
	\mA_{k-1}:\left\lvert
	\begin{array}{cccc}
		&\mathbb{R}^{\Fq} \times\ldots \times\mathbb{R}^{\Fq} & \longrightarrow & \mathbb{R}^{\Fq} \\
		& (\Mu_0,\ldots,\Mu_{k-1}) &\mapsto& \mA_{k-1}(\Mu_0,\ldots,\Mu_{k-1}).
	\end{array}
	\right.
\end{equation}
and the error is governed by the mesh size $h$ and degree $p$:
\begin{equation}
	\label{formula1}
	\U_{k} = \Mu_{k}+\delta^h \Mu_{k},\quad \text{such that} \quad \|\delta^h \Mu_{k}\| \approx O_{k}(h^{\max (2,p-1)}). 
\end{equation}
Though, $\delta^h \Mu_{k}$ refers to the error related only to the discrete form of the equation if all $\Mu_0,\ldots,\Mu_{k-1}$ are exactly known.

Let us consider that we have approximations $\U_0, \U_1,\dots \U_{k-1}$ of terms $\Mu_0$, $\Mu_1,\dots,\Mu_{k-1}$, respectively, up to order related to $h$ and $p$, \emph{i.e.} following \eqref{formula1}.
For simplicity, we denote by $\MU_{k-1} := (\Mu_0, \ldots,\Mu_{k-1})^\T$ the vector of the first $k$ exact terms of the series, and by $\MU^h_{k-1} := (\U_0,\ldots,\U_{k-1})^\T$ its discrete form in the triangulation of the given \ac{FEMF}. Replacing $\MU_{k-1}$ by $\MU^h_{k-1}$ in the right-hand side of \eqref{sys_exct_un} and denoting by $\bU^h_{k}$ its solution:
\begin{eqnarray*}
	\M^h \bU^h_{k} &=& \frac{1}{k}\mA_{k-1}(\MU^h_{k-1}).
\end{eqnarray*}
Using $\MU^h_{k-1}=\MU_{k-1}+ \delta^h \MU_{k-1}$ and the perturbation theory, we can write:
\begin{eqnarray*}
	\M^h \bU^h_{k}& = \displaystyle\frac{1}{k}\mA_{k-1}(\MU_{k-1} ) &+~~ \frac{1}{k}\bD \mA_{k-1}(\MU_{k-1} ).\delta^h \MU_{k-1} + O\Big((\delta^h \MU_{k-1})^2\Big)\\
	& =  a_{k-1} &+~~ \delta^h \mA_{k-1}+  O\Big((\delta^h \MU_{k-1})^2\Big)\\
	& =  a_{k-1} &+~~ \delta^h a_{k-1},
\end{eqnarray*}
where $\bD \mA_{k-1}(\MU_{k-1} )$ is a liner operator and is the Fréchet differential of $\mA_{k-1}$, assuming it is differentiable. Though, the unknown $\bU^h_k$ could be estimated using $\U_k$ and an error $\delta^h\U_k$ related to the perturbation $\delta^h a_{k-1}$ if $\MU_{k-1}^h$ is used:
\begin{equation*}
	\bU^h_k = \U_k  + \delta ^h\U_k,\quad \delta^h \U_k = (\M^h)^{-1} \delta^h a_{k-1}.
\end{equation*}
If $\mA_{k-1}$ is muti-linear operator, therefore $ O\Big((\delta^h \MU_{k-1})^2\Big)$ is negligible. If $\|\delta^h \MU_{k-1}\|\ll 1$, though $(\delta^h \MU_{k-1})^2$ could be neglected in front of $\delta^h \mA_{k-1}$. In both cases, we have according to the condition number $\kappa(\M^h)$:
\begin{eqnarray*}
	\frac{\|\delta^h\U_k\|}{\|\U_k\|}
	&\leqslant & {\kappa(\M^h)} \cdot\frac{\|\delta^h a_{k-1}\|}{\|a_{k-1}\|}\\
	&\leqslant &   \frac{\kappa(\M^h)}{k}\cdot \frac{\op{\bD \mA_{k-1}(\MU_{k-1} )}\cdot\,\|\delta^h \MU_{k-1}\|}{\|a_{k-1}\|}.
\end{eqnarray*}
We replace $\|\U_k\|$ by $\|(\M^h)^{-1}a_{k-1}\|$, we can state:
\begin{eqnarray}
	\label{rel_delta_uk}
	{\|\delta^h\U_k\|}\leqslant \frac{\kappa(\M^h)}{k} \cdot\op{(\M^h)^{-1}}
	\cdot\op{\bD \mA_{k-1}(\MU_{k-1} )}\cdot\|\delta^h \MU_{k-1}\|.
\end{eqnarray}
This means that $\kappa(\M^h)\cdot\op{\bD \mA_{k-1}(\MU_{k-1} )}\cdot\op{(\M^h)^{-1}}$ is the amplification factor of  $\|\delta^h \MU_{k-1}\|$.
We recall that $ \bU^h_k = \Mu_k + \delta^h \Mu_k  + \delta^h \U_k$, where $\Mu_k$ is the exact solution, $\delta^h \Mu_k$ is the error related to the precision of the FEM in case where all $\Mu_0, \ldots,\Mu_{k-1}$ are known, and $\delta^h \U_k$ is the amplification of $\delta^h \MU_{k-1}$ by the condition number of the mass matrix $\M^h$.
We compute now $\bU_{k+1}$ using $\U_0,\ldots,\U_{k-1}$ and $\bU^h_k$:
\begin{eqnarray*}
	\M^h \bU^h_{k+1} &=& \frac{1}{k+1}\mA_{k}(\U_0,\ldots,\U_{k-1},\bU^h_k)\\
	(\text{perturbing $\mA_{k-1}$})&=&\frac{1}{k+1} \left(\mA_{k}( \MU_k)+
	\bD \mA_k(\MU_k).\delta^h \MU_k
	+ \frac{\partial \mA_{k}}{\partial \Mu_k}\cdot\delta^h \U_k+...\right)\\
	&=& a_{k} + \delta^h \mA_k + \frac{1}{k+1}  \frac{\partial \mA_{k}}{\partial \Mu_k}\cdot\delta^h \U_k + O\big((\delta^h\MU_k)^2+(\delta^h \U_k)^2 \big).
\end{eqnarray*}
The second equality results by, after perturbing $\mA_k$, the fact that $\cfrac{\partial \mA_{k}}{\partial \Mu_k}(\Mu_0,..,\Mu_k)$ is a multi-linear operator.
Considering:
\begin{equation*}
	\bU^h_{k+1} = \Mu_{k+1} + \delta^h \Mu_{k+1}+\delta^h \U_{k+1}+\delta^h \bU^h_{k+1},
\end{equation*}
where $\delta^h \Mu_{k+1}$ is the error relative to the discrete form of the operators and having used exact solutions of all precedents $\Mu_k$, $\delta^h \U_{k+1}$ and $\delta^h \bU^h_{k+1}$ are the errors amplification of $\delta^h\MU_{k}$ and $\delta^h\U_k$ respectively, such that:
\begin{eqnarray*}
	\frac{\|\delta^h \U_{k+1}\|}{\|\U_{k+1}\|} &\leqslant & \frac{\kappa(\M^h)} {k+1}\cdot\op{\bD \mA_{k}(\MU_k)}\cdot\frac{\| \delta^h \MU_{k}\|}{\|a_k\|},\\
	\frac{\|\delta^h \bU^h_{k+1} \|}{\|\U_{k+1}\|} &\leqslant & \frac{\kappa(\M^h)}{k+1}\cdot\op{\frac{\partial \mA_{k}}{\partial \Mu_k}(\MU_k)}\cdot\frac{ \|\delta^h \U_k\|}{\|a_{k} \|}.
\end{eqnarray*}
Two cases are presented here:
\begin{itemize}
	\item If $\|\delta^h\U_k\| < \|\delta^h \Mu_k \|$, thus the propagation of $\delta^h \U_k$ to $\bU^h_{k+1}$, denoted by $\delta^h\bU^h_{k+1}$:
	\begin{equation*}
		{\|\delta^h \bU^h_{k+1} \|}\leqslant \frac{\kappa(\M^h)}{k+1}\cdot \op{(\M^h)^{-1}} \cdot \op{\frac{\partial \mA_{k}}{\partial \Mu_k}(\MU_k)}  \cdot{ \|\delta^h \U_k\|},
	\end{equation*}
	is neglected compared to that of $\delta^h \Mu_k$, denoted by $\delta^h\U_{k+1}$:
	\begin{eqnarray*}
		{\|\delta^h\U_{k+1}\|} &\leqslant & \frac{\kappa(\M^h)}{k+1}\cdot\op{(\M^h)^{-1}} \cdot {\op{\bD \mA_{k}(\MU_k)}\cdot\| \delta^h \MU_{k}\|},
	\end{eqnarray*}
	which is stable if the following is verified for every $k$:
	\begin{eqnarray}
		\label{rel_kappa_Mh_Ak}
		\frac{\kappa(\M^h)}{k}  \cdot \op{(\M^h)^{-1}} \cdot\op{\bD \mA_{k-1}(\MU_{k-1} )} < &1, & \forall k.
	\end{eqnarray}

	\item If $\|\delta^h\U_k\| \geqslant \|\delta^h \Mu_k \|$,we can conclude that:
	\begin{eqnarray*}
		\frac{\kappa(\M^h)}{k} \cdot \op{\bD \mA_{k-1}(\MU_{k-1} )} &\geqslant &\frac{\|\delta^h\Mu_k\|}{\| \delta^h \MU_{k-1}\|}\cdot\frac{\|a_{k-1}\|}{\|\delta^h\U_k\|},
	\end{eqnarray*}
	with $\|a_{k-1}\|\geqslant \op{(\M^h)^{-1}}^{-1} \cdot\|\U_k\|$ and ${\|\delta^h\Mu_k\|}\sim{\| \delta^h \MU_{k-1}\|}$ 	though,
	\begin{equation*}
		\frac{\kappa(\M^h)}{k}  \cdot \op{\bD \mA_{k-1}(\MU_{k-1} )} \cdot\op{(\M^h)^{-1}} \geqslant 1,
	\end{equation*}
	and the error is represented by:
	\begin{eqnarray*}
		\|\delta^h \bU^h_{k+1} \| &\leqslant & \frac{ \kappa(\M^h)}{k+1}\cdot
		\op{(\M^h)^{-1}} \cdot \op{\frac{\partial \mA_{k}}{\partial \Mu_k}(\MU_k)}\cdot \|\delta^h \U_k\|.
	\end{eqnarray*}
	Replacing $\|\delta^h \U_k\|$ by its formula \eqref{rel_delta_uk}, we have the error $\delta^h\MU_{k-1}$ is amplified to $\bU^h_{k+1}$ by the factor
\begin{equation*}
	\frac{ \kappa(\M^h)^2}{k\times(k+1)}\cdot\op{(\M^h)^{-1}}^2 \cdot
	\op{\frac{\partial \mA_{k}}{\partial \Mu_k}}\cdot\op{\bD \mA_{k-1}(\MU_{k-1} )}.
\end{equation*}

\end{itemize}

In the next Section, we present the technique that we apply to reduce these errors by modifying a given linear system to reduce its condition number.

\nomenclature{$\delta^h\U_k$}{Propagation of $\delta^h\MU_{k-1}$ to $\U_k$}
\nomenclature{$\delta^h\bU_{k+1}$}{Propagation of $\delta^h\U_{k}$ to $\bU_{k+1}$}
\nomenclature{$\bU_{k+1}$}{Approximation of $\Mu_{k+1}$ using $\U_0,\ldots,\U_k$}

\section{The idea}
\begin{figure}[!h]
\centerline{
	\includegraphics[width=0.7\textwidth]{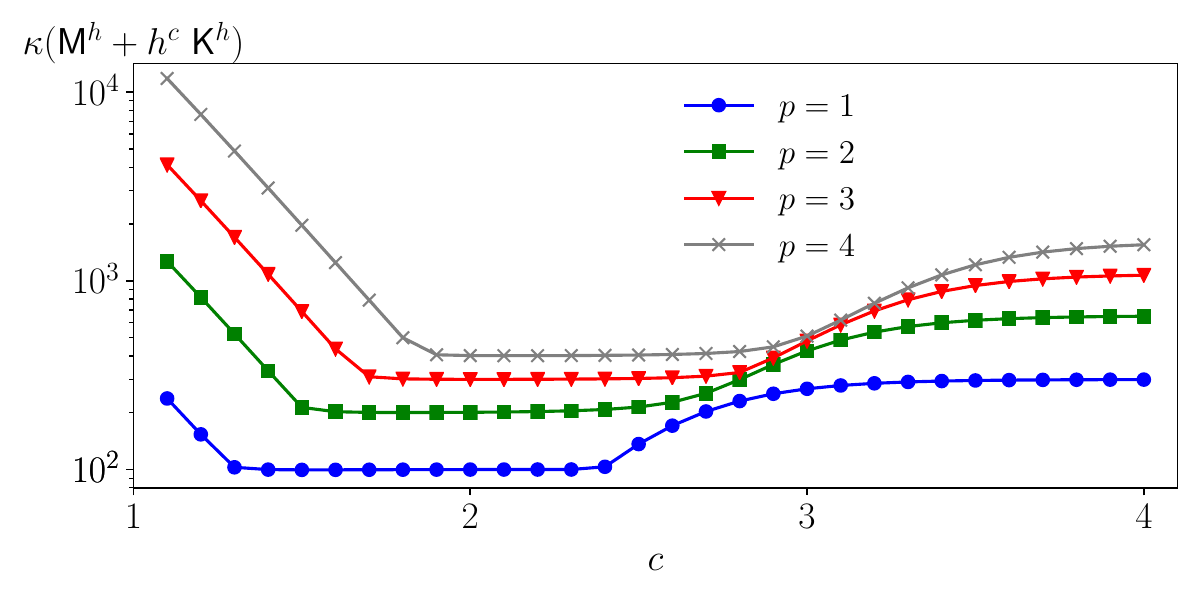}
}
\caption{Optimality of the functional in problem \eqref{prob_argmin_alpha} relative to $\alpha=h^c$ for $N_x=100$}
\label{fig61}
\end{figure}

After showing that the condition number of the mass matrix is the major source of the amplification error in computing $\U_k$, our idea consists of modifying the linear system to reduce its condition number, thus to reduce amplification of errors along $\U_k$.
First, we will modify the linear system and solve a new one as follows:
\begin{equation}
\label{modified_ls}
(\M + \alpha \K) \Mu = b,
\end{equation}
where $\K$ is a given matrix, which will be taken as the stiff matrix for the moment (we explain its advantage later on), and $\alpha$ is a scalar parameter that we consider it depends on the mesh size $h$ and degree $p$ of the \ac{FEMF}. The choice of $\alpha$ will be considered as it may minimize the condition number of the modified matrix; \emph{i.e.} for a give size mesh $h$ and a degree $p$, find $\alpha$ that is solution of the following minimization problem:
\begin{equation}
\label{prob_argmin_alpha}
\argmin\limits_{\alpha} \kappa \big(\M^h + \alpha \K^h\big).
\end{equation}
To show the optimality, we present in Figure \ref{fig61} the evolution of the condition number $\kappa \big(\M^h + \alpha \K^h\big)$ relative to the one dimensional geometry when $N_x=100$ and $p\in \{1,2,3,4\}$ and for a range of values of $\alpha$. We can see in this figure that a global minimum could be found.
We present also in Table \ref{tab3} values of $c$ corresponding to $\alpha = h^c$ that minimizes the condition number $\kappa\big(\M^h + \alpha \K^h\big)$ for different \ac{FEMF} associated to the geometry $\Omega=[0,1]$.
\begin{table}[!h]
\caption{Values of $c$ relative to optimal choice $\alpha = h^c$.}
\label{tab3}
\begin{center}
	\begin{tabular*}{\textwidth}{@{\extracolsep\fill}lllllll@{}}
		\toprule
		& $h=1/50$ & $h=1/100$ & $h=1/200$ & $h=1/300$ & $h=1/400$ & $h=1/500$ \\
		\midrule
		$p=1$ & 1.96 & 1.96 & 1.96 & 1.96 & 1.96 & 1.96\\
		$p=2$ & 2.18 & 2.02 & 1.9  & 1.86 & 1.82 & 1.80 \\
		$p=3$ & 2.52 & 2.28 & 2.14 & 2.06 & 2.02 & 2.00 \\
		$p=4$ & 2.68 & 2.4  & 2.22 & 2.14 & 2.10  & 2.06\\
		\bottomrule
	\end{tabular*}
\end{center}
\end{table}
According to it, the optimal choice of $\alpha \sim h^{2 \pm c_1(h,p)}$: for small size of mesh, $\alpha$ decreases when $p$ increases, and slightly decreases when $h$ decreases. It seems to approaching the value $h^2$; \emph{i.e.} $c_1(h,p)\xrightarrow[]{ h\rightarrow 0}0$. In a considered case with a fixed $h$ and $p$, what we have to do is, once and for all, solving the problem \eqref{prob_argmin_alpha} and store $\alpha$ for further computation.

\begin{remark}
In the case of computing terms of series in the \ac{FEMF} through the recurrence Formula \eqref{strong-form}, adding $\alpha \K $ in the left hand side of the linear system \eqref{lin_sys_heat} is equivalent of having a new recurrence formula:
\begin{equation}
	\begin{array}{llcc}
		u_{k+1} - \alpha_k \Delta u_{k+1} &= \displaystyle\frac{\nu}{k+1}  \Delta u_k& \Longrightarrow&  (\M +\alpha_k \K)\U_{k+1} =a_{k}.
	\end{array}
\end{equation}
Though, the new variational formulation associated to the new approach will be as follows:
\begin{equation}
	\int_{\Omega} u_{k+1}v\,dx + \alpha_k\int_{\Omega} \nabla u_{k+1}\nabla v\,dx = -\frac{\nu}{k+1}\int_{\Omega} \nabla u_k \nabla v\,dx, \quad \forall k>0.
\end{equation}
\end{remark}
Before testing the effect of this technique on lowering instabilities in computing terms of the series, we consider that the initial linear system is the discrete form of the minimization problem:
\begin{equation}
\min\limits_{v \in V} \mathcal{J}^k(v) \Longleftrightarrow  \M\Mu = b.
\end{equation}
Now we add the following quantity
\begin{equation*}
E^k_{\alpha_k}(v_{}) =\frac{\alpha_k}{2}  \int_{\Omega} |\nabla v|^2\,dx
\end{equation*}
to the functional $\mathcal{J}^k$. It is equivalent to ask the solution to be with minimal total variation. Thus, the modified linear system will be the discrete form of a new minimization problem:
\begin{equation*}
\begin{array}{c}
	\min\limits_{v \in V} \mathcal{J}^k_{\alpha_k}(v) \Longleftrightarrow (\M+\alpha_k \K)\Mu = b , \quad \text{such that}\quad
	\displaystyle \mathcal{J}^k_{\alpha_k}(v) :=\mathcal{J}^k(v) +E^k_{\alpha_k}(v_{}).
\end{array}
\end{equation*}
This is why, the matrix $\K$ is chosen in the modified system. The new solution is of minimal total variation.

In the proposed technique, the coefficient $\alpha_k$ is considered depending on the mesh size $h$ and the degree $p$. We will check if there is any possible dependence of $\alpha_k$ on $k$.

\subsection{Optimal choice of $\alpha_k$}
\begin{figure}[!h]
	\centerline{
		\includegraphics[width=0.45\textwidth]{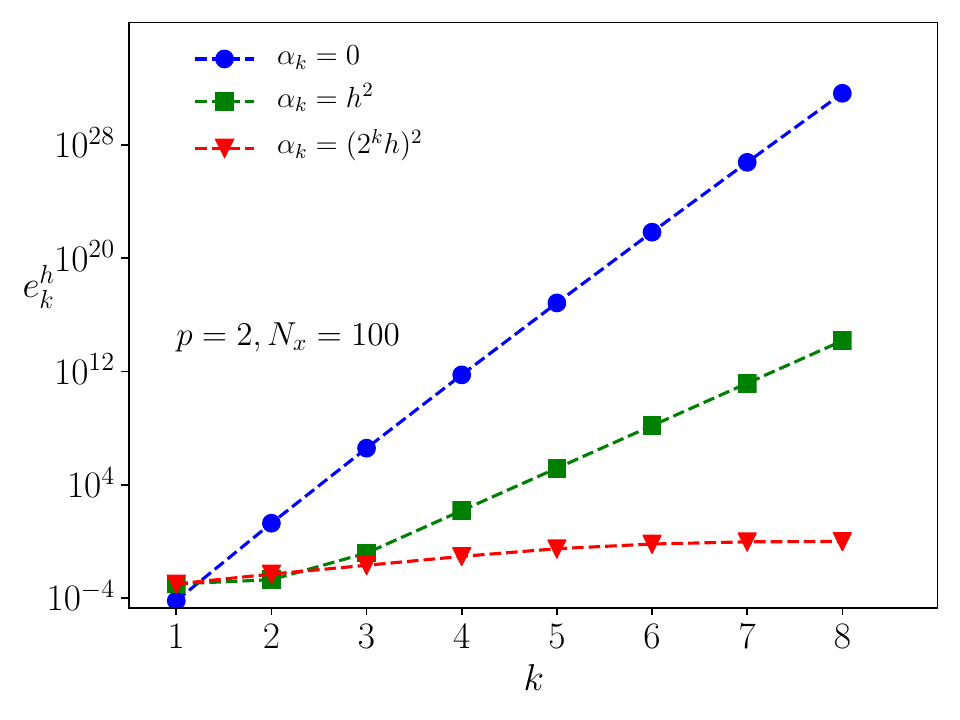}
		\includegraphics[width=0.45\textwidth]{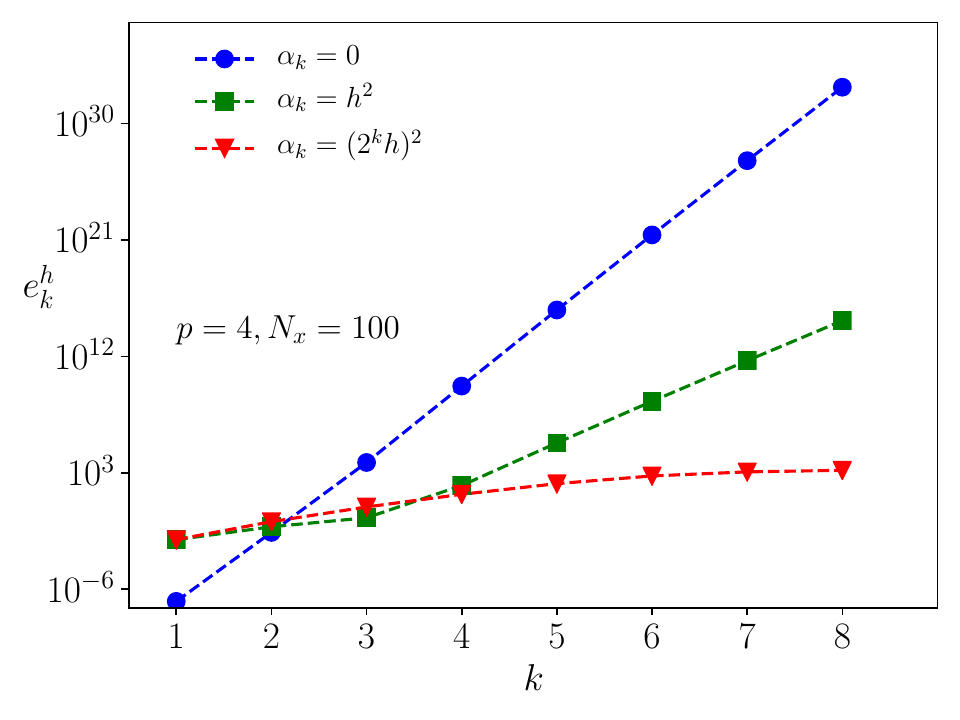}
	}
	\caption{Errors $\Me_k^{h,2}$ when applying stabilization techniques with different patterns of $\alpha_k$.}
	\label{fig6}
\end{figure}

For the moment, we do not have any idea on how to choose $\alpha_k$ when $k$ changes; \emph{i.e.} how to choose $\alpha_k$ when advancing in computing terms of the series. The optimal choice would be the one which minimizes the following objective function:
\begin{equation}
	\label{prob_argmin_alpha1}
	\begin{array}{c}
		\argmin\limits_{\alpha_k} \displaystyle \int_{\Omega} |u^h_{k+1} - u_{k+1}|^2\\
		\text{subject to}\\
		(\M^h+ \alpha_k \K^h) \U_{k+1} = \frac{1}{k+1} \mA_k(\U_0,\ldots,\U_k).
	\end{array}
\end{equation}
In reality, we do not have the exact solution, thus it is not possible to use this function. First of all, we will assume that $\alpha_k = h^c$ does not depend on $k$ and will take the value from Table \ref{tab3} for every $h$ and $p$. Second, we assume that $\alpha_k \approx (2^k h)^c$; \emph{i.e.} the coefficient of the artificial diffusion we add here will be doubled when advancing from $u_k$ to $u_{k+1}$. Following these two assumptions, we plot in Figure \ref{fig6} the error between exact and numerical terms series obtained using those two patterns of $\alpha_k$. We will add the one when $\alpha_k = 0$ as a reference case.
It turns out that the pattern $\alpha_k = (2^k h)^2$ presents a better approximation for terms $u_k$, $\forall\,k\geqslant 2$ than the constant pattern. This confirms that for some choice of $\alpha_k$ depending on $k$, it is more stable when $\alpha_k$ is constant. But how to find the optimal choice? Can we find an expression for the ratio ${\alpha_{k}}/{\alpha_{k-1}}$ ?

We consider the case when we are at the level of computing $\bU^h_k$ using the stabilization technique with $\alpha_{k-1}$:
\begin{equation*}
	(\M^h + \alpha_{k-1}\K^h)\bU^h_k = \frac{1}{k} \mA_k(\U_0,\ldots,\U_{k-1}) ,
\end{equation*}
where we have $\bU^h_k = \U_k + \delta^h\U_k + \delta^h \U_k$, such that $\delta^h\U_k$ is related to perturbation term $\delta^h a_{k-1}$ and is neglected in front of $\delta^h \U_k$. We have
\begin{equation}
	\label{magnif_factor_k}
	\|\delta^h\U_k\| \leqslant \frac{\kappa_{k-1}}{k}  \left[ \frac{\|\U_k\|}{\|a_{k-1}\|}\right] \cdot \op{ \bD \mA_{k-1}(\MU_k)}\cdot\|\delta^h\MU_{k-1}\|,
\end{equation}
with $\kappa_{k-1} :=\kappa(\M^h+\alpha_{k-1} \K^h)$. Then, with a new stabilization coefficient $\alpha_k$, we will approximate $\Mu_{k+1}$ by $\bU^h_{k+1}$ as follows:
\begin{eqnarray*}
	(\M^h + \alpha_{k}\K^h)\bU^h_{k+1} &= & \frac{1}{k+1} \mA_{k+1}(\U_0,\ldots,\U_{k-1},\bU^h_{k}).
\end{eqnarray*}
We follow the same technique as in the last section and we estimate the amplification of the error $\delta^h\Mu^h_k$ in $\bU^h_{k+1}$, denoted by $\delta^h\bU^h_{k+1}$. We can find that, for any recurrence formula related to $\mA_k$, the magnification factor of the perturbation $\delta^h\Mu^h_k$ reads:
\begin{equation}
	\label{minimizing_factor}
	\frac{\kappa_{k-1}\kappa_{k}}{k(k+1)}\cdot \op{(\M^h + \alpha_{k-1}\K^h)^{-1}} \cdot \op{(\M^h + \alpha_{k}\K^h)^{-1}}.
\end{equation}
Hence, minimizing this factor will minimize the amplification of the error to the next iteration. This means that for a given $\alpha_{k-1}$ used to compute $\Mu^h_k$, the coefficient $\alpha_k$ is chosen such that it minimizes problem \eqref{minimizing_factor}.

\nomenclature{$\mathcal{J}^k$}{Functional relative to $\Mu_k$}
\nomenclature{$\alpha_k$}{$k$\up{th} stabilization coefficient}

\subsection{Link with the \acf{DMP}}
The maximum principle of the one-dimensional Heat equation, $u_t-\nu u_{xx}=0$ defined on $[0,1]\times [0,T]$, states that, if the source term is null, the solution has its maximum only in the initial condition, at the left $\{(t,x), x=0\, \&\, t\in [0,T]\}$ or at the right boundaries $\{(t,x), x=1\, \& \,t\in [0,T]\}$. The \ac{DMP} \cite{DMP_Varga,ciarlet_DMP,DMP_Farago,chen_DMP} states that the numerical approximation of the discrete probelm should satisfy the maximum principle of the problem, \emph{i.e.} the numerical solution should admits its maximum either in the initial condition or at the boundaries.
Consider the discrete problem approximating the solution over the grid decomposed into inner points $x^{\inn}$ and boundaries $x^\rb$. We denote by $u^\inn(t)\coloneqq u(t,x^\inn)$ and $u^\rb(t)\coloneqq u(t,x^\rb)$ the solution at the inner and boundaries points respectively, and $\Mbu(t)=\big(\Fbu^1(t),\dots,\Fbu^\Fq(t) \big)^{\top}$ the vector containing functions in time that approximate the solution at every inner $x^\inn$ or boundary node $x^\rb$ of the mesh. Then the matrix-vector form of the discrete heat equation in the framework of FEM can be written in the following form:
\begin{equation}
 \M\frac{d}{dt}\Mbu(t) = \nu \K \Mbu \Longleftrightarrow
 \left[
 \begin{array}{cc}
\M^\inn & \M^\rb \\0&I
 \end{array}
 \right]
 \frac{d}{dt}\left[
 \begin{array}{c}
\Mbu^\inn \\ \Mbu^\rb \end{array}
 \right] =
 \left[
 \begin{array}{cc}
\nu \K^\inn & \nu \K^\rb \\0&I
 \end{array}
 \right]
 \left[
 \begin{array}{c}
\Mbu^\inn \\ \Mbu^\rb
 \end{array}
 \right],
\end{equation}
where $\M^\inn$ and $\M^\rb$ are matrix block containing elements of the mass matrix $\M$ and $\K^\inn$ and $\K^\rb$ are relative to the stiff matrix $\K$. The computation of the terms of the series of the discrete problem consists of solving the following linear problem many times until reaching a defined rank $m$:
\begin{equation}
 \M\Mu_{k+1} = -\frac{\nu}{k+1} \K\Mu_k \quad \forall k\in \{0,...,m\} \Longleftrightarrow \Mu_k = \frac{(-\nu)^k}{k!} \left(\M^{-1}\cdot \K\right)^k \Mu_0.
\end{equation}
When the partial sum is applied as numerical integrator of the problem, the solution is approximated as follows:
\begin{equation}
 \Mbu(t) \approx \sum\limits_{k=0}^m \Mu_kt^k =  \left[ \sum\limits_{k=0}^m\frac{(-\nu)^kt^k}{k!} \left(\M^{-1} \cdot \K\right)^k  \right] \Mu_0.
\end{equation}
We use the form of the inverse of $\M$ to write the following:
\begin{equation}
\sum\limits_{k=0}^m \Mu_kt^k =  \left[ \sum\limits_{k=0}^m\frac{(-\nu)^kt^k}{k!} \left(
 \left[
 \begin{array}{cc}
(\M^\inn)^{-1} & - (\M^\inn)^{-1}\cdot\M^\rb \\0&I
 \end{array}
 \right]
 \left[
 \begin{array}{cc}
 \K^\inn &  \K^\rb \\0&I
 \end{array}
 \right]
\right)^k  \right]
\left[
\begin{array}{c}
 \Mu_0^\inn\\\Mu_0^\rb
\end{array}
\right],
\end{equation}
thus we multiply both matrices to have:
\begin{equation}
 \sum\limits_{k=0}^m \Mu_kt^k =  \left[
 \sum\limits_{k=0}^m\frac{(-\nu)^kt^k}{k!} \left(
  \begin{array}{cc}
   (\M^\inn)^{-1}\cdot\K^\inn &(\M^\inn)^{-1}\cdot(\K^\rb-\M^\rb)\\
   0 &I
  \end{array}
 \right)^k
 \right] \left[
\begin{array}{c}
 \Mu_0^\inn\\\Mu_0^\rb
\end{array}
\right].
\end{equation}
We can show that:
\begin{equation}
  \left(\M^{-1} \K\right)^k =
 \left(
  \begin{array}{cc}
   \big((\M^\inn)^{-1}\cdot\K^\inn\big)^k & \sum\limits_{p=0}^{k-1} \big((\M^\inn)^{-1}\cdot\K^\inn\big)^p\cdot(\M^\inn)^{-1}\cdot(\K^\rb-\M^\rb)\\
   0 &I
  \end{array}
 \right),
\end{equation}
then the solution is approximated by the partial sum as follows:
\begin{equation}
 \Mbu(t) \approx \left[
 \begin{array}{cc}
\sum\limits_{k=0}^m\frac{(-\nu)^kt^k}{k!}\big((\M^\inn)^{-1}\cdot\K^\inn\big)^k & \sum\limits_{k=0}^m\frac{(-\nu)^kt^k}{k!}\sum\limits_{p=0}^{k-1} \big((\M^\inn)^{-1}\cdot\K^\inn\big)^p\cdot(\M^\inn)^{-1}\cdot(\K^\rb-\M^\rb) \\
0& \sum\limits_{k=0}^m\frac{(-\nu)^kt^k}{k!}
 \end{array}
 \right] \left[
\begin{array}{c}
 \Mu_0^\inn\\\Mu_0^\rb
\end{array}
\right]
\end{equation}
If the problem is considered with homogenous Dirichlet boundary condition, the approximation of the solution $\Mbu^\inn$ at the inner part is simplified by the following:
\begin{equation}
 \Mbu^\inn(t) \approx \sum\limits_{k=0}^m\frac{(-\nu)^kt^k}{k!}\big((\M^\inn)^{-1}\cdot\K^\inn\big)^k \Mbu_0^\inn,
\end{equation}
where the maximum principle in this case states that the initial condition is the function that presents the maximum of the solution over all the time; \emph{i.e}:
\begin{equation}
 \left\|\sum\limits_{k=0}^m\frac{(-\nu)^kt^k}{k!}\big((\M^\inn)^{-1}\cdot\K^\inn\big)^k \Mbu_0^\inn \right\|\leqslant \left\|\Mbu_0^\inn \right\|.
\end{equation}
Therefore, we can preserve the maximum principle of the solution obtained via the partial sum of the discrete problem by verifying the following inequality:
\begin{equation}
  \left\|\sum\limits_{k=0}^m\frac{(-\nu)^kt^k}{k!}\big((\M^\inn)^{-1}\cdot\K^\inn\big)^k  \right\|\leqslant 1.
\end{equation}
This mean that the smaller the condition number of the mass matrix $\M$ is the smaller the norm of the above operator of summation is, as $\kappa(\M) = \left\|\M\right\|\cdot\|\M^{-1}\|$. Our proposition consider of replacing $\M^\inn$ by $\M^\inn + \alpha_k\K^\inn$, and choosing $\alpha_k$ that minimize the condition number of $\kappa(\M^\inn + \alpha_k\K^\inn)$, allowing a bigger value of $t$.
In conclusion, the proposed stabilization helps satisfying the \ac{DMP} of the heat equation.

\subsection{Empirical simplification proposition: $\alpha_{k+1} = R(h,p)\alpha_k$}
\label{sec_prop}
It is costly to find $\alpha_k$ at every iteration of computing $\U_k$ and for every problem that minimizes the magnification factor established in formula \eqref{minimizing_factor}. Therefore, we consider here that the ratio $\alpha_{k+1}/\alpha_k$ will depend only on the \ac{FEMF}; \emph{i.e.} the size of mesh $h$ and degree $p$:
\begin{equation}
	\frac{\alpha_{k+1}}{\alpha_k} \simeq R(h,p).
\end{equation}
Finding this constant will be obtained using formula \eqref{minimizing_factor}. To proceed, we first compute $\alpha_0$ while minimizing the magnification factor of $\delta^h\Mu^h_k$ in formula \eqref{magnif_factor_k}; \emph{i.e.}
\begin{equation}
	\label{find_alpha0}
	\alpha_0 = \argmin\limits_{\alpha}\kappa(\M^h + \alpha\K^h).
\end{equation}
Then we consider that $\alpha_k = R(h,p) \alpha_{k-1}$ and find $R(h,p)$ that minimizes the magnification factor in formula \eqref{minimizing_factor}, \emph{i.e.}
\begin{equation}
	\label{find_Chp}
	R(h,p) := \argmin\limits_{r}\kappa(\M^h + r\alpha_0\K^h)\, \op{\left( \M^h + r\alpha_0\K^h\right)^{-1}}.
\end{equation}
\begin{figure}[!h]
	\centerline{
		\includegraphics[width=0.45\textwidth]{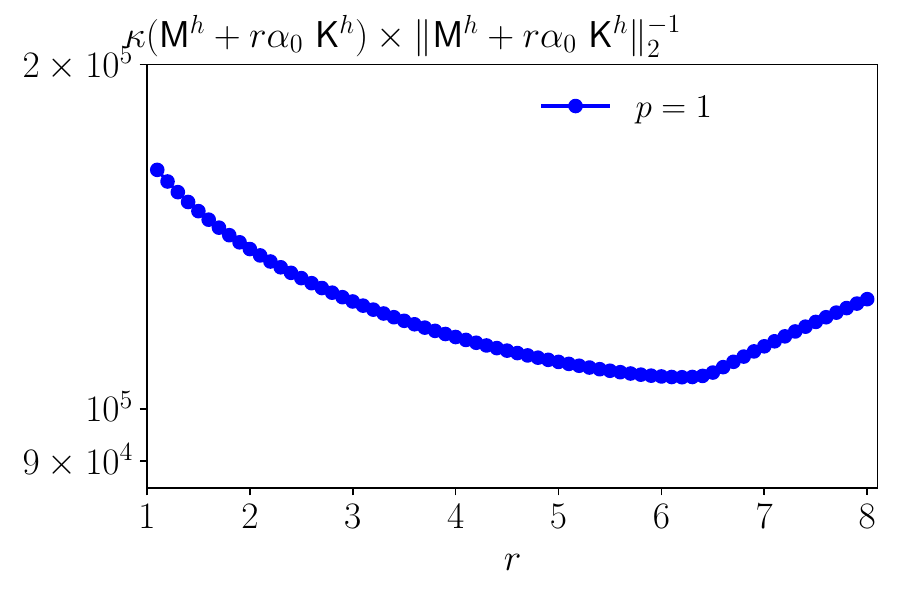}
		\includegraphics[width=0.45\textwidth]{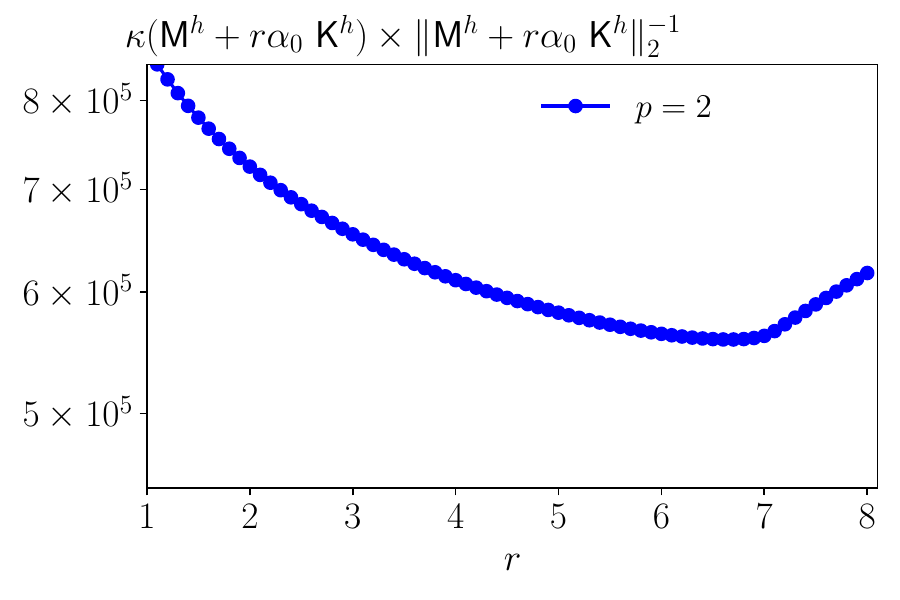}
	}
	\centerline{
		\includegraphics[width=0.45\textwidth]{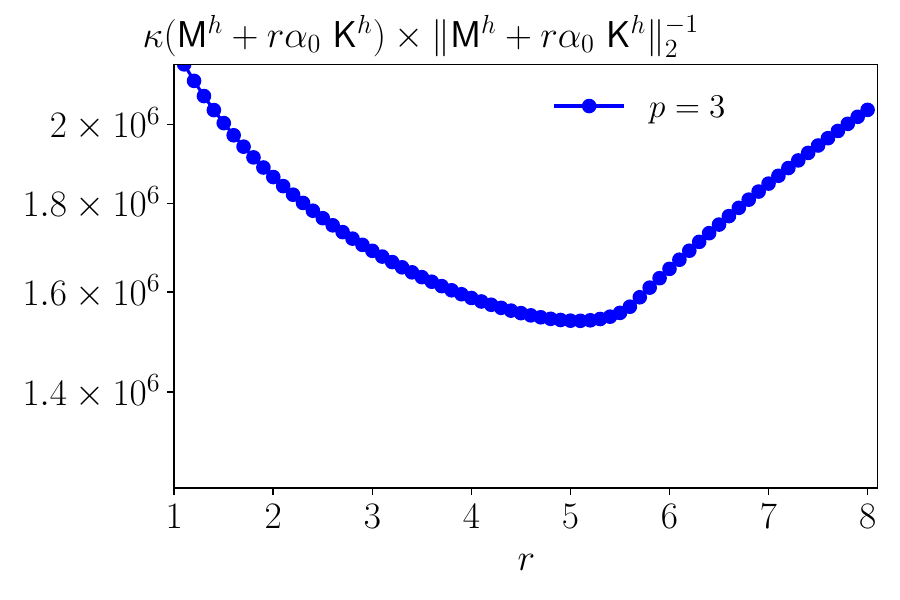}
		\includegraphics[width=0.45\textwidth]{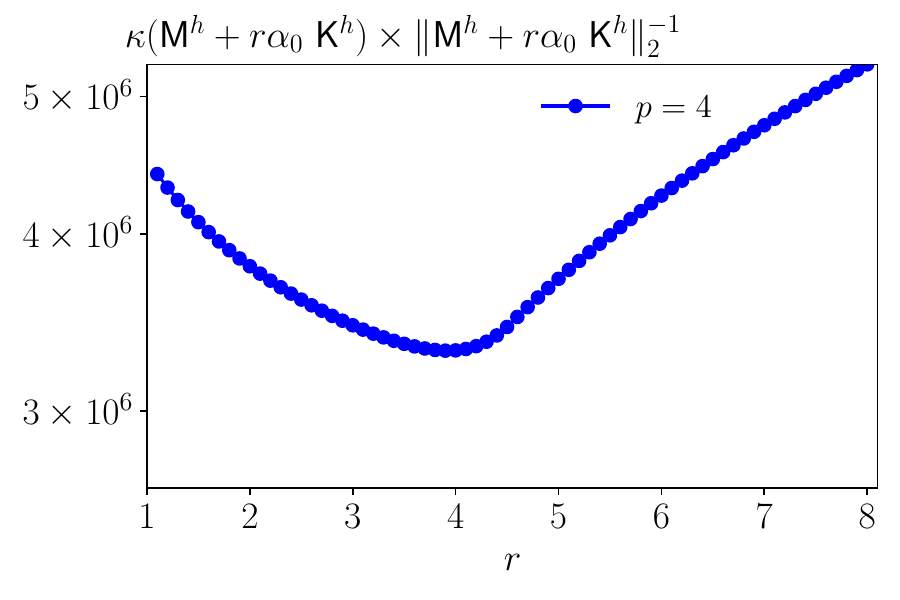}
	}
	\caption{Optimality of the functional in problem \eqref{find_Chp} relative to the space formulation associated tp the geometry $[0,1]$}
	\label{fig62}
\end{figure}
In Figure \ref{fig62}, we present the evolution of the functional relative to $r$ for a test case when $N_x = 100$ and for different degree $p=\in\,\{1,2,3,4\}$, where the minimum is the optimal choice $R(h,p)$. These minima are presented later in Table \ref{tab-heat-alpha-1d} for the one dimensional \ac{FEMF}.

\subsection{Algorithm for Computing $\Mu^h_k$}
\label{sec:alg}
We present below the steps to follow in order to approximate $u_k$ by $\U_k$ that verify the recurrence Formula \eqref{recurrence_formula_pde}. For every \ac{FEMF}, \emph{i.e.} for every mesh discretization with size $h$ and an associated approximation space formulation $\mathbb{V}_{h,p}$ with degree $p$, do:
\begin{enumerate}
	\item Establish, for every $k$, the form of $\mA_k$ relative to $\A_k$.
	\item Construct the mass matrix $\M^h$ and stiff matrix $\K^h$ defined in \eqref{mass_stiff_matrix}.
	\item Find $\alpha_0$ and $R(h,p)$ as in \eqref{find_alpha0} and \eqref{find_Chp}, and consider $\alpha_k = R(h,p)\alpha_{k-1}$.
	\item Having $\U_0$, compute $\U_k$ for $k=1,\ldots,m$ by solving the above linear systems:
	\begin{equation}
		(\M^h + \alpha_k\K^h)\U_{k+1} = \frac{1}{k+1}\mA_k(\U_0,\ldots,\U_k)
	\end{equation}
\end{enumerate}
Once these terms are computed, we can apply the \ac{BPL} algorithm from Diagram \ref{alg:borel_pade_laplace} for $\sum\limits_{k=0}^{m}\U_k t^k$ and approximate the solution for $u(t,x^i), x^i \in \mathcal{T}_h$. To advance in time, we follow Algorithm \ref{alg:1} in Section \ref{sec_cont}.

\section{Numerical experiments}

In this part, we proceed to apply the \ac{BPL} algorithm presented in Section \ref{sec_cont}, after computing the series terms using algorithm in Section \ref{sec:alg}, on integrating in time some celebrated examples of \ac{PDE}s in the \ac{FEMF}. We limit our examples on one dimensional geometry. The examples we are using to test the presented procedure are:
\begin{enumerate}
	\item Heat equation with homogeneous Dirichlet boundary conditions.
	\item Burgers equations.
\end{enumerate}

\subsection{Some notations}
We present a \ac{PDE} by its semi-discrete formulation:
\begin{equation*}
	\M^h \frac{\d \,\Mu(t)}{\d t} = \mA \big(\Mu(t)\big), \quad \Mu_0:=\Mu(0).
\end{equation*}
For an approximation $\U(t)$ of the above IVP, we denote by $\Res(t)(\Mu_0)$ the residual error defined as follows:
\begin{equation*}
	\Res(t)(\Mu_0) = \left\| \M^h \frac{\d \U(t)}{\d t} - \mA \big(\U(t)\big) \right\|\cdot \left\|\frac{\d\,\U(t)}{\d t}\right\|^{-1}.
\end{equation*}
This residual error will be evaluated to rate the efficiency of the approximation as we consider that we do not have the exact solution. It is also used in the \ac{BPL} algorithm to estimate the largest time step $\Delta t_n^{\BPL}$ in the continuation process of \ac{BPL} algorithm. To this end, we denote $\RES$ the global residual error over all the segment $[0,T]$:
\begin{equation}
	\label{Res}
	\RES(t) = \int_0^{t}\Res(\tau)(\Mu_0)\, \d \tau.
\end{equation}

\subsection{Heat equation}
We take again the equation in Section \ref{sec_case-study} with initial condition $u_0(x) = \sin(\pi x)$ where its time series solution has $u_k(x) = \frac{(-\nu \pi^2)^k}{k!}\sin(\pi x)$. This is a convergent series. However, partial sums could not be efficient as $\|u_k\|$ increase first before decreasing after a defined rank related to $\nu$.

In the \ac{FEMF} and for fixed value of $N_x =50,100,200$, we construct several spaces $\mathbb{V}_{h,p}$ for $p=1,2,3,4$ to approximate $u_k(x)$ in the associated FEM formulation. After constructing the mass matrix and the stiff matrix, we write the discrete form of the recurrence formula that compute vectors $\U_k$:
\begin{equation*}
	\M^h \U_k =- \frac{\nu}{k} \K^h \U_{k-1}.
\end{equation*}
In this case, $\mA_k = \K^h$ is a linear operator and the error $\delta^h \U_k$ can be bounded as
\begin{equation*}
	\|\delta^h \U_k\|\leqslant \frac{\nu}{k}\cdot\kappa(\M^h) \cdot\op{\M^h}^{-1}\cdot\op{\K^h\U_{k-1}} \cdot \|\delta^h \Mu_{k-1}\|,
\end{equation*}
for what the propagation error to the $k$-th term $\U_k$ is stable if and only if:
\begin{equation}
	\label{stab_cond_uk_heat}
	\kappa(\M^h) \op{\M^h}^{-1}\cdot\op{\K^h} <\frac{k}{\nu}.
\end{equation}
The smaller the diffusion coefficient is, the more we can increase the size of the mesh or the degree $p$ while having a stable computation. We present in Table \ref{tab6} values of the logarithm in base ten of $\kappa(\M^h) \cdot\HS{\M^h}^{-1}\cdot\HS{\K^h}$ for different values of $h$ and $p$.
\begin{table}[!h]
	\caption{Values of $\log_{10}\left(\kappa(\M^h) \HS{\M^h}^{-1}\cdot\HS{\K^h}\right)$ for different \ac{FEMF}.}
	\label{tab6}
	\begin{center}
		\begin{tabular*}{\textwidth}{@{\extracolsep\fill}lllllllll@{}}
			\toprule
			  & $h=1/20$ & $h=1/50$ & $h=1/100$ & $h=1/200$ \\ \midrule
			$p=1$ & 3.35& 5.10& 5.77& 6.39\\
			$p=2$& 4.99& 5.81& 6.43& 7.03\\
			$p=3$& 5.40& 6.21& 6.81& 7.42\\
			$p=4$& 5.70& 6.53& 7.14& 7.75\\
			\bottomrule
		\end{tabular*}
	\end{center}
\end{table}
These values are upper bounds of those of slopes $\slope(h,p)$ relative to the growth rate of errors $\Me_k^{h,p}$ in the terms of series of the heat equation presented in Table \ref{tab1}. In addition, we remark that for $p\in\,\{2,3,4\}$ ,the logarithmic rate increases by $2$ when $h$ is divided by $100$, which is also present in the slopes in Figure \ref{fig1}.
We proceed to modify the linear system and add the artificial diffusion $\alpha_k \K^h$ in the left hand side by following the procedure presented in Section \ref{sec_prop}.
We present in Table \ref{tab-heat-alpha-1d} values of $\alpha_0$ and $R(h,p)$ relative to each formulation.
\begin{table}[!h]
	\caption{Results of values $c$ and $R(h,p)$ for $\alpha_0 = h^c$ and $\alpha_{k+1} = R(h,p)\alpha_k$ for different \ac{FEMF} associated to $\Omega=[0,1]$.}
	\label{tab-heat-alpha-1d}
	\begin{center}
		\begin{tabular*}{\textwidth}{@{\extracolsep\fill}lllllllll@{}}
			\toprule
			& \multicolumn{2}{@{}c}{$h=1/20$} & \multicolumn{2}{@{}c}{$h=1/50$} & \multicolumn{2}{@{}c}{$h=1/100$} & \multicolumn{2}{@{}c}{$h=1/200$}\\
			\cmidrule{2-3}\cmidrule{4-5}\cmidrule{6-7}\cmidrule{8-9}
			&  $c$ & $R(h,p)$ & $c$ & $R(h,p)$ & $c$ & $R(h,p)$ & $c$ & $R(h,p)$   \\
			\midrule
			$p=2$ &  2.0 & 1.9 & 1.86 & 2.85 & 1.8 & 3.8 & 1.76 & 5.4 \\
			$p=3$ &  2.3 & 2.1 & 2.12 & 3.5 & 2.0 & 4.7 & 1.94 & 6.3 \\
			$p=4$ &  2.5 & 1.7 & 2.26 & 2.7 & 2.1 & 3.4 & 2.02 & 4.4\\
			$p=5$ & 2.8 & 1.9 & 2.4 & 2.1 & 2.3 & 3.5  &  2.1& 3 \\
			\bottomrule
		\end{tabular*}
	\end{center}
\end{table}
For every \ac{FEMF} (for a given $h$ and $p$), we compute both values: $c$ and $R(h,p)$ in order to modify the system we solve for computing $\U_k$.

\begin{figure}[!h]
	\centerline{
		\includegraphics[scale=0.45]{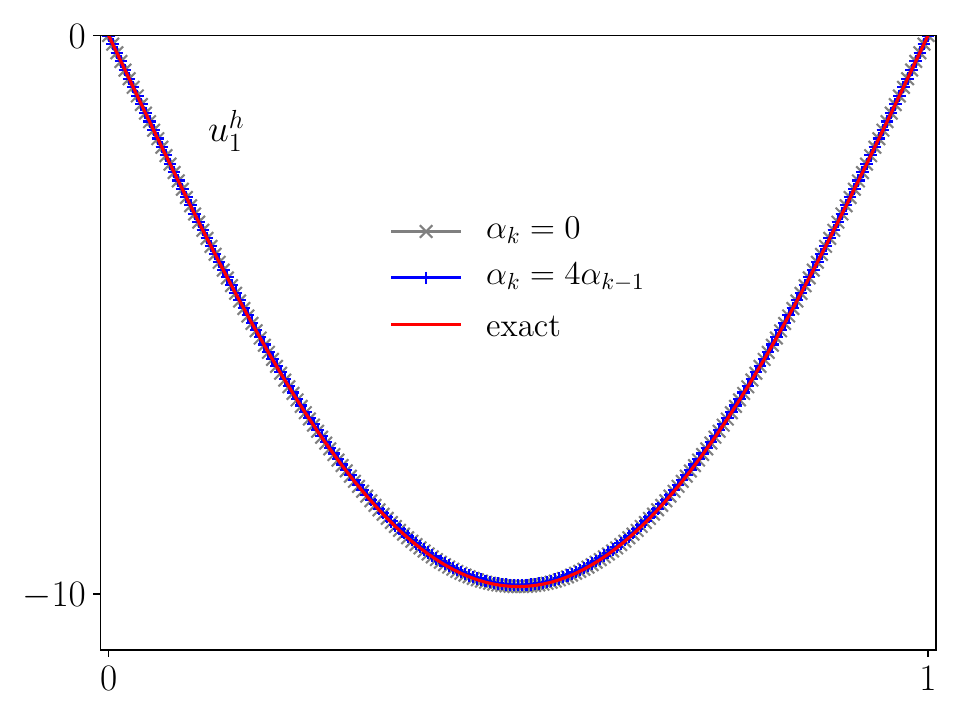}
		\includegraphics[scale=0.45]{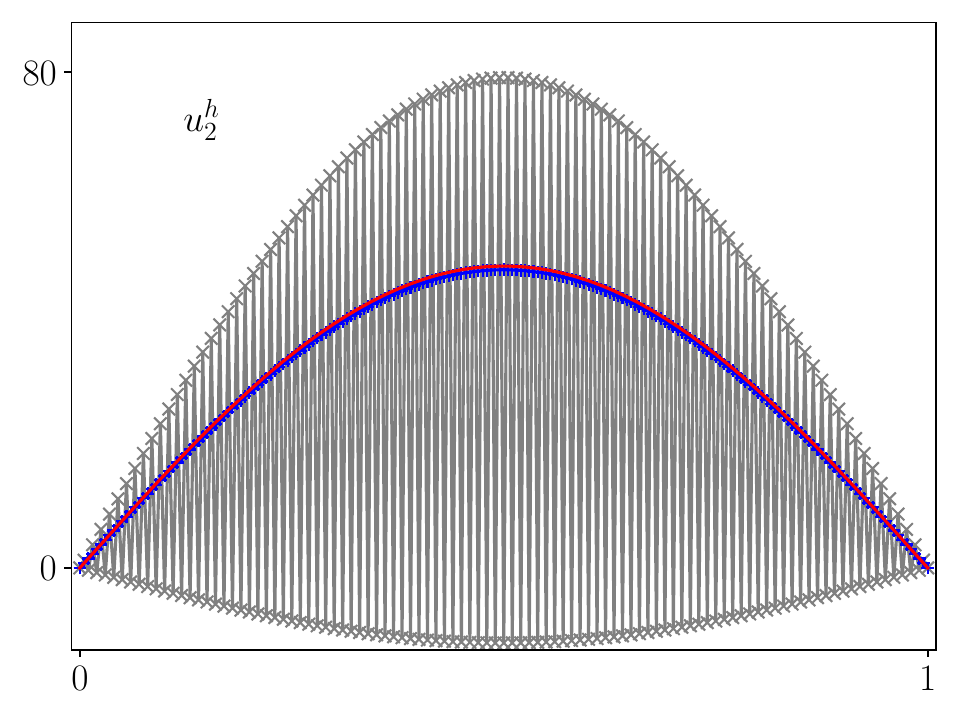}
	}
	\centerline{
		\includegraphics[scale=0.45]{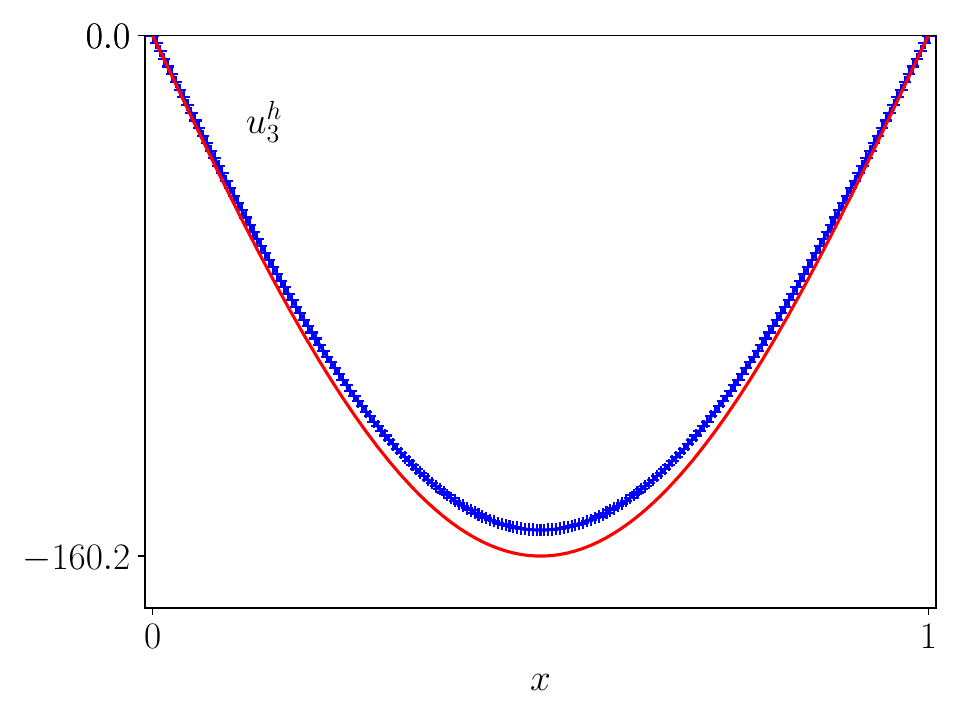}
		\includegraphics[scale=0.45]{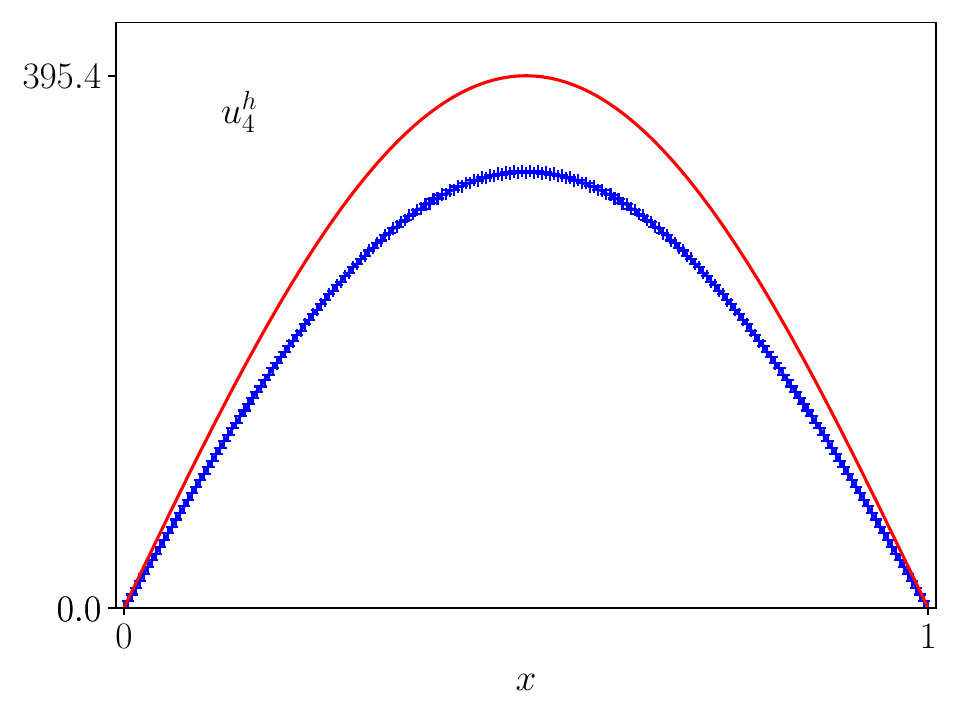}
	}
	\caption{Plots of the exact terms $\Mu_k, \,k\in\{ 1,2,3,4\}$ of the heat equation's series, their approximations in the space $\mathbb{V}_{h,p}$ for $h= 1/100$ and $p = 2$, with and without stabilization.}
	\label{fig7}
\end{figure}
We plot in Figure \ref{fig7} approximations of terms $u_1,u_2,u_3,u_4$ obtained in the framework of $h=1/100$ and $p=2$, using the proposed process of stabilisation. We compare them with the exact terms and the ones obtained without the stabilisation ($\alpha_k=0$). We can see that approximation of $u_1$ is valid with or without the stabilization. However, approximation of $ u_2$ without the proposed technique oscillates around the exact solution. This oscillation will lead to the explosion in computing the third terms as it is not shown in this figure. Nevertheless, the stabilization technique produces a stable approximation of $u_2$ such that we are able to compute the third and fourth terms as shown in Figure \ref{fig7}.

\begin{figure}[!h]
	\centerline{
		\includegraphics[scale=0.45]{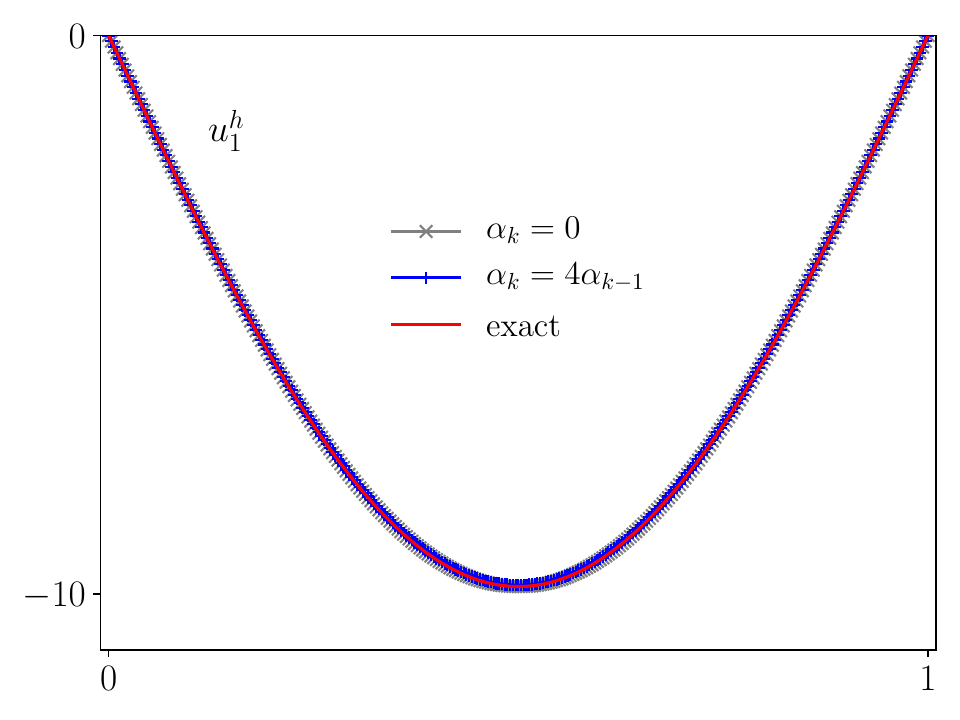}
		\includegraphics[scale=0.45]{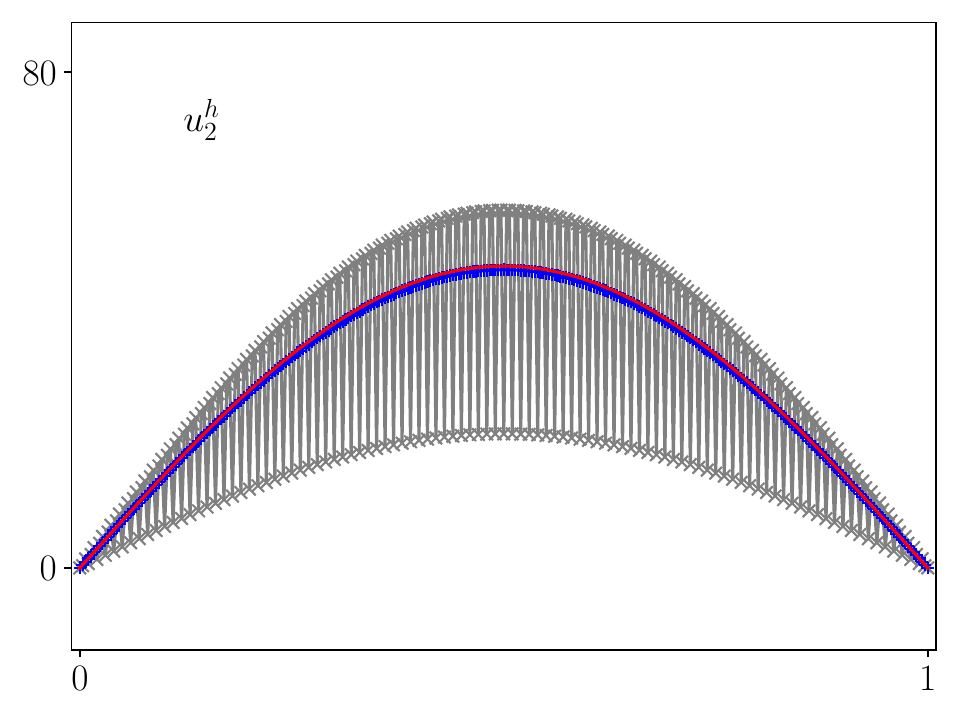}
	}
	\caption{Plots of the exact terms $\Mu_k, \,k\in\{ 1,2\}$ of the heat equation's series, their approximations in the space $\mathbb{V}_{h,p}$ for $h= 1/100$ and $p = 3$, with and without stabilization.}
	\label{fig8}
\end{figure}
\begin{figure}[!h]
	\centerline{
		\includegraphics[scale=0.45]{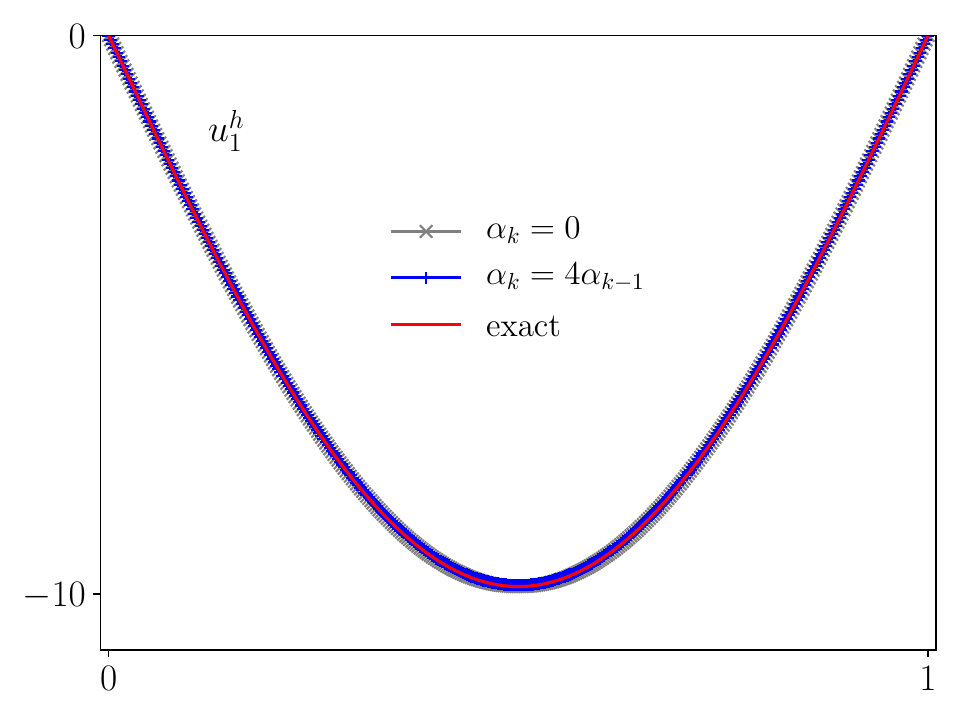}
		\includegraphics[scale=0.45]{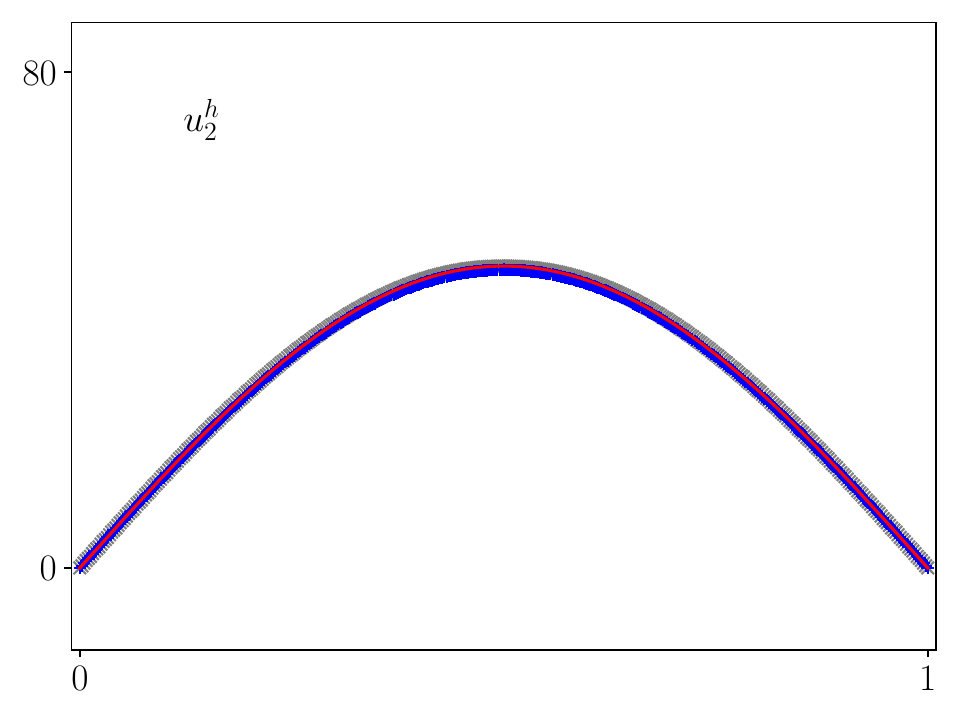}
	}
	\centerline{
		\includegraphics[scale=0.45]{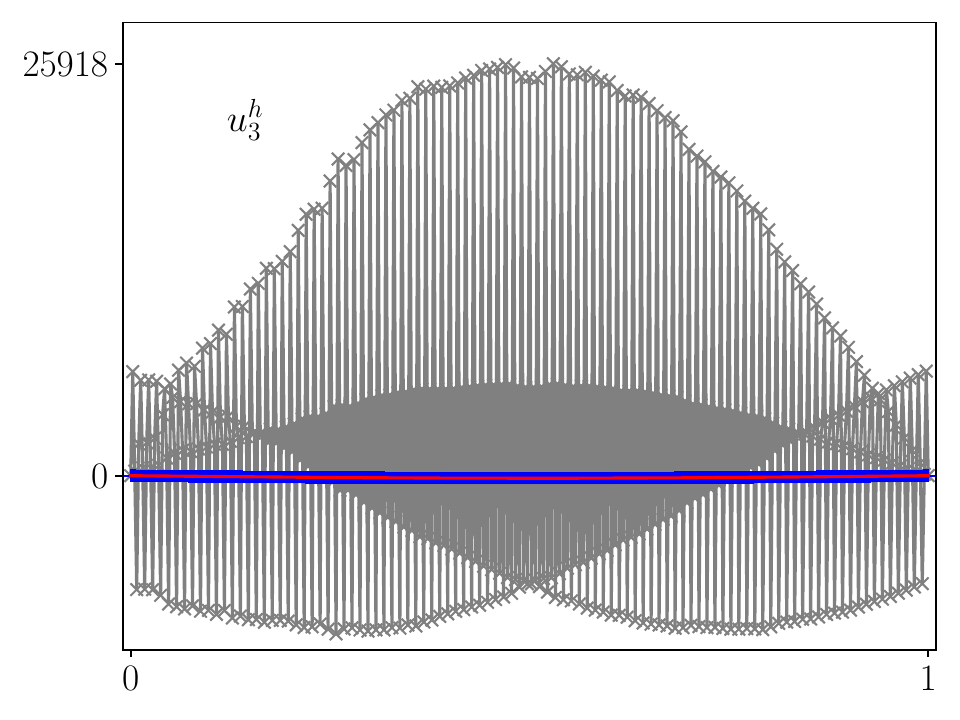}
		\includegraphics[scale=0.45]{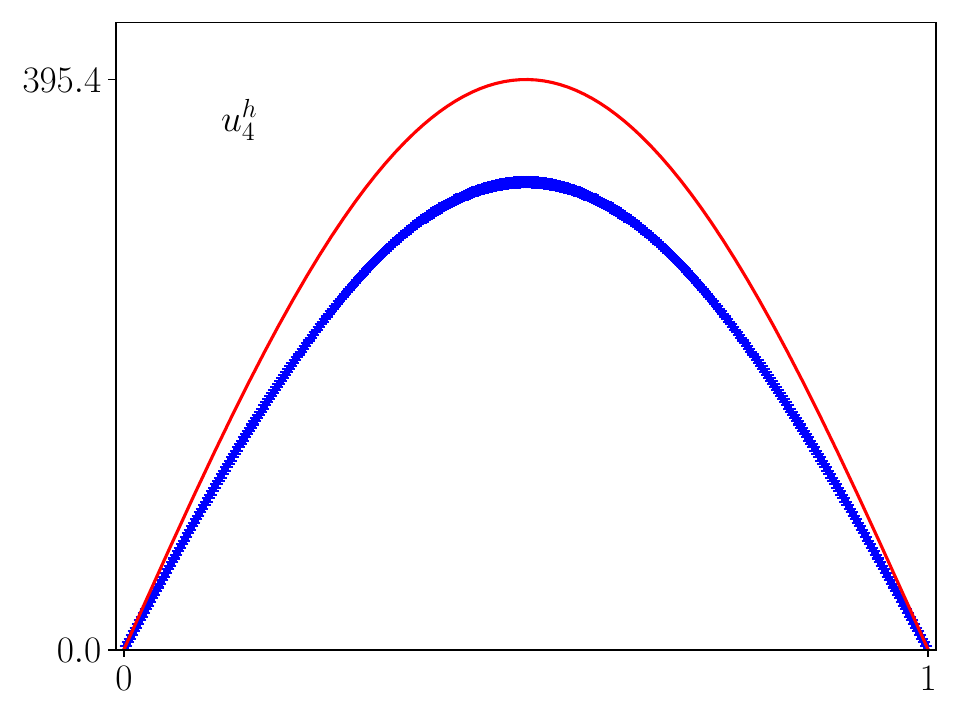}
	}
	\caption{Plots of the exact terms $\Mu_k, \,k\in\{ 1,2,3,4\}$ of the heat equation's series, their approximations in the space $\mathbb{V}_{h,p}$ for $h= 1/100$ and $p = 4$, with and without stabilization.}
	\label{fig9}
\end{figure}

We have computed these terms, in the same triangulation but with $p=3$ with and without the stabilization process. We present these approximations in Figure \ref{fig8} for the case whether or not the stabilization process is used. We can remark that the oscillation in the second term persists for $p=3$ (Figure \ref{fig8}) but its magnitude is less than $p=2$. It seems that, for $p=4$ (Figure \ref{fig9}) and whether or not the stabilization is applied, the approximation of the second term is valid, but the one for the third approximation is not valid when stabilization is not applied.

After showing how the proposed stabilization process is beneficial for having a stable approximations of terms $\U_k$, we will proceed to apply the \ac{BPL} algorithm.

We consider nine \ac{FEMF}s that come form combination of $N_x\in\,\{ 20,50,100\}$ with the degree of elements $p\in\,\{1,2,3\}$ and we run the simulation for a fixed time step that is presented below with all others parameters needed to run the \ac{BPL} algorithm:
\begin{align*}
	m &  =5 ,& r=s&=2,&  \Delta t^{\BPL} &= 5\times 10^{-3}, &N_g & =20, & T & =1.
\end{align*}
The residual error $\Res(\U)(t_n)$ at points $t_n = n\Delta t^{\BPL}$ is computed, and $\RES$ is approximated by a trapezoidal method. Results are reported in Table \ref{tab4}.
We can see that the global residual error decreases as $h$ decreases or $p$ increases.
\begin{table}[!h]
	\caption{Evaluation of the residual $\RES(1)$ from Formula \eqref{Res} for different configurations of the \ac{FEMF} for the one dimensional heat equation.}
	\label{tab4}
	\begin{center}
		\begin{tabular*}{\textwidth}{@{\extracolsep\fill}llll@{}}
			\toprule
		    & $h=1/20$ & $h=1/50$ & $h=1/100$ \\ \midrule
			$p=1$ & 0.0573 &  0.0357 & 0.0251 \\
			$p=2$ & 0.0363 & 0.0225 &  0.0159\\
			$p=3$ & 0.0274  &  0.0170 & 0.0120 \\
			\bottomrule
		\end{tabular*}
	\end{center}
\end{table}

We conclude the effectiveness of using a stabilisation technique to compute terms of the series in the case of heat equation, and allow us to use resummation techniques as the \ac{BPL} algorithm. Next, the Burgers equation is investigated too in the following section.

\subsection{Burger's equation}

In this section, we show how the proposed approach to solve numerically the one-dimensional Burgers equation given can be applied:
\begin{equation*}
	\partial_t u + u\, \partial_x (u)=\nu \,\partial_x^2 u, \quad (t,x) \in (0,T)\times ]0,1[.
\end{equation*}
After writing the solution as a formal power series in time, the recurrence formula associated to this \ac{PDE} is given by:
\begin{equation*}
	u_{k+1} = \frac{1}{k+1} \left(\nu\,\partial_x^2 u_k-\sum\limits_{i=0}^k u_i\,\partial_x( u_{k-i})  \right),
\end{equation*}
where the corresponding weak formulation is shown as follows:
\begin{equation*}
	a(u_{k+1}^h,v) =   \frac{-1}{k+1}\left( \nu\ell_1(u_k,v) +\sum\limits_{i=0}^{k}\ell_2(u^h_i,u^h_{k-i},v)\right),
\end{equation*}
with
\begin{align*}
	a(u,v) & =  \int_{\Omega} u \,v\,\d x, &\ell_1(u,v) &= \int_{\Omega} \partial_x u \,\partial_x v \,\d \,x,   & \ell_2(u,w,v) &= \int_{\Omega} u\,\partial_x(w) \,v \,\d x.
\end{align*}

In the \ac{FEMF}, the linear system associated to the unknown $\U_{k+1}$ is given below:
\begin{equation*}
	\M^h \U_{k+1} = \frac{1}{k+1}\mA_k(\U_0,\ldots,\U_{k}) =
	\frac{1}{k+1}\left( -\nu \K^h \U_k - \sum\limits_{r=0}^{k} \U_r:( \D^h:\U_{k-r}) \right),
\end{equation*}
where $\D^h$ is a tensor of rank three defined by:
\begin{equation*}
	\D_{i,j,l}^h =\int_{\Omega}\phi^i\, \phi^j_x\, \phi^l \,\d x.
\end{equation*}
Thus, the amplification factor in computing the Burgers' series terms, $\U_{k}$, is given below:
\begin{equation}
	\label{ampl_fact_burg}
	\frac{\kappa(\M^h)}{k} \cdot \HS{\M^h}^{-1}\cdot \HS{ \frac{\partial \mA_{k-1}}{\partial{\Mu_k}}} \leqslant
	\frac{\kappa(\M^h)}{k} \cdot\HS{\M^h}^{-1}\cdot\left( \nu \HS{\K^h}+ k \HS{\D^h}\right).
\end{equation}

\subsection{Inviscid Burger's equation}
We consider the case when $\nu=0$. For the initial condition $u_0(x) = \sin(2\pi x)$, we give below the formal solution of $u_k, k\in \{1,\ldots,6\}$
\begin{eqnarray*}
	u_1(x)& = & -\pi \sin(4 \pi x),\\
	u_2(x)& = & \frac{\pi^2}{2} \left( 3 \sin(6\pi x) - \sin(2\pi x) \right),   \\
	u_3(x)& = & -\frac{\pi^3}{3!} \left( 16 \sin(8\pi x) - 8\sin(4\pi x) \right),\\
	u_4(x)& = &  \frac{\pi^4}{4!} \left( 125 \sin(10\pi x) - 81 \sin(6\pi x) +2 \sin(2\pi x) \right),                \\
	u_5(x)& = & -\frac{\pi^5}{5!} \left( 1296 \sin(12\pi x) - 1024 \sin(8\pi x) + 80 \sin(4\pi x) \right), \\
	u_6(x)& = & \frac{\pi^6}{6!} \left( 16807 \sin(14\pi x) - 15625 \sin(10\pi x) + 2187\sin(6\pi x) - 5 \sin(2\pi x) \right) .
\end{eqnarray*}
\begin{table}[htp]
	\caption{Values of $\log_{10}\Big(\kappa(\M^h) \cdot\HS{\M^h}^{-1}\cdot\HS{\D^h}\Big)$ for different \ac{FEMF}}
	\label{tab7}
	\begin{center}
		\begin{tabular*}{\textwidth}{@{\extracolsep\fill}lllll@{}}
			\toprule
			 & $h=1/20$ & $h=1/50$ & $h=1/100$ & $h=1/200$ \\ \midrule
			$p=1$& 0.99 & 2.13& 2.35& 2.52 \\
			$p=2$& 2.20& 2.42& 2.58& 2.73 \\
			$p=3$& 2.34& 2.55& 2.70& 2.85 \\
			$p=4$& 2.45& 2.68& 2.83& 2.99 \\
			\bottomrule
		\end{tabular*}
	\end{center}
\end{table}

\begin{figure}[!h]
	\centerline{
		\includegraphics[scale=0.45]{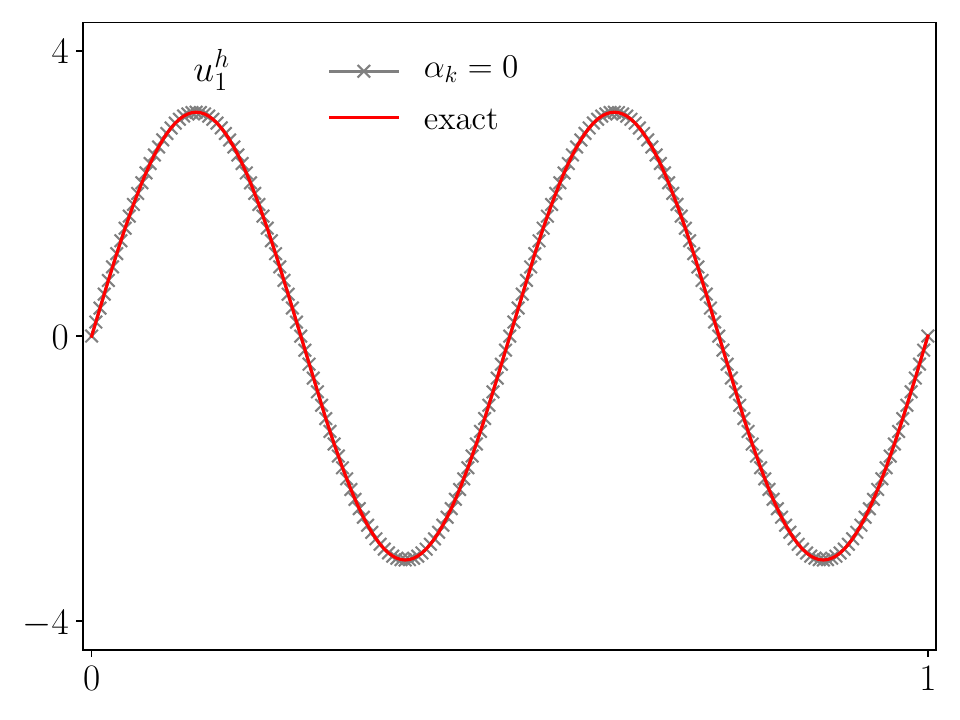}
		\includegraphics[scale=0.45]{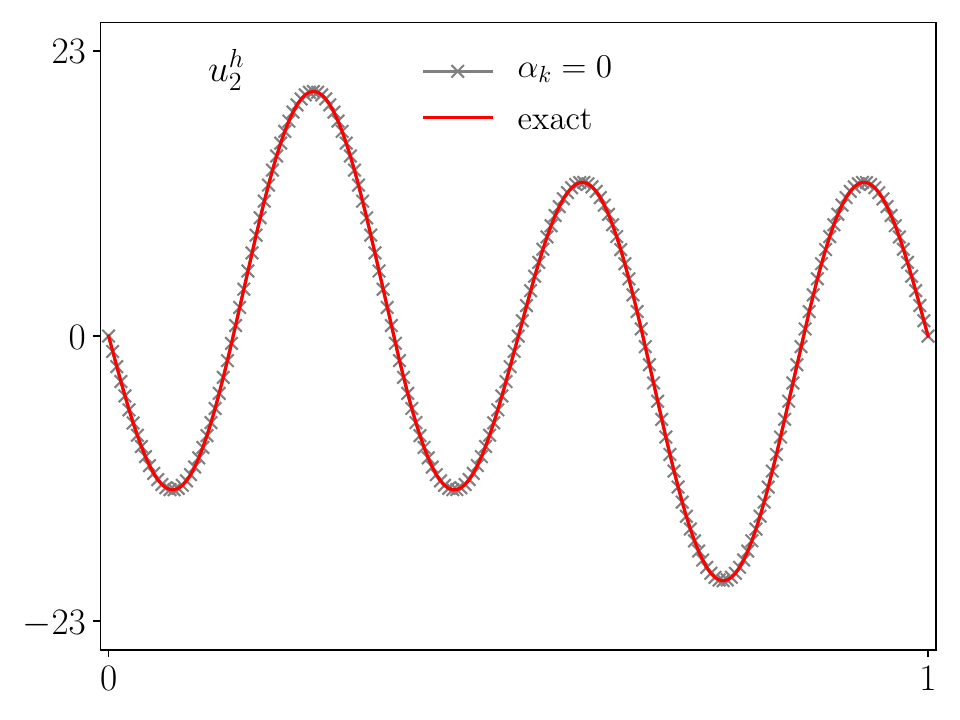}
	}
	\centerline{
		\includegraphics[scale=0.45]{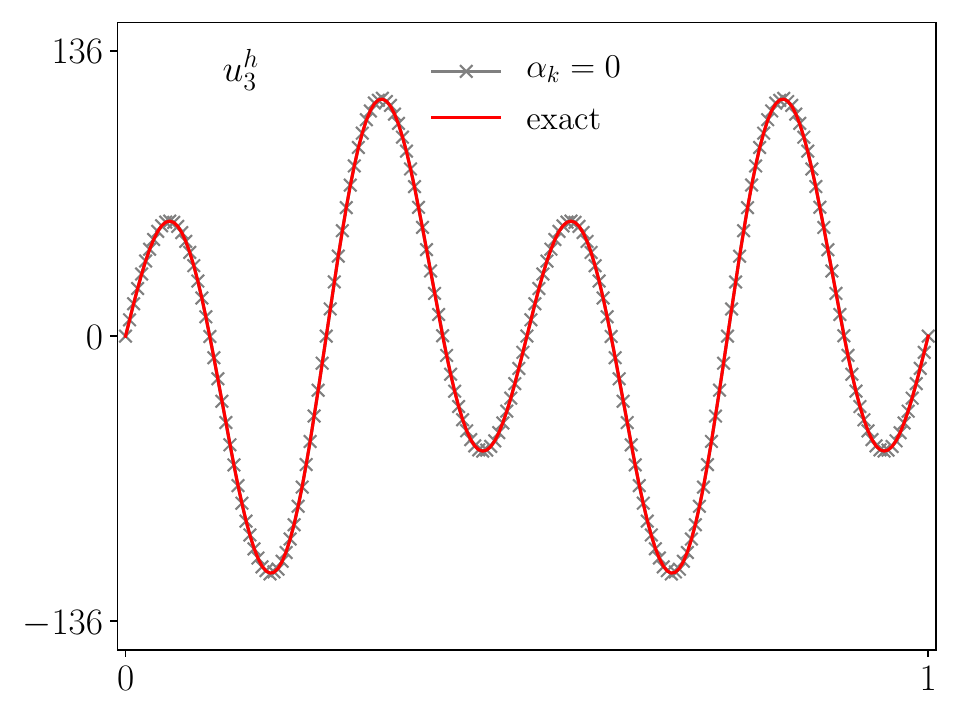}
		\includegraphics[scale=0.45]{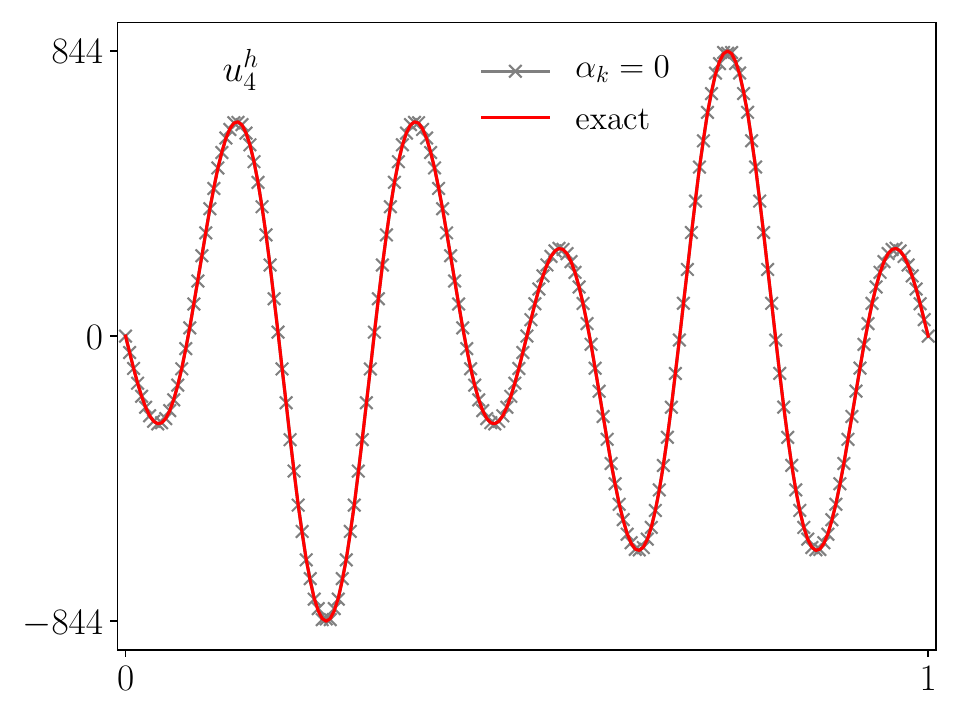}
	}
	\caption{Plots of the exact terms $\Mu_k, \,k\in\{ 1,2,3,4\}$ of the Burger equation's series, their approximations in the space $\mathbb{V}_{h,p}$ for $h= 1/100$ and $p = 2$, without stabilization.}
	\label{fig10}
\end{figure}

In this case, the amplification factor is given by $\kappa(\M^h) \cdot \HS{\M^h}^{-1}\cdot\HS{\D^h}$, where we present in Table \ref{tab7} its logarithmic values in base ten. We remark that these values are much smaller then those, presented in Table \ref{tab6}, and relative to the magnification factor of the heat equation.  We remark also that they do not increase as fast as they do in Table \ref{tab6} when $h$ decreases or $p$ increases.
In this case we compute approximations of terms $\Mu_k$ without the artificial diffusion term $\alpha_k \Delta \Mu_k$ and we plot them in Figure \ref{fig10} when $p=2$ and $h=1/100$.

\subsection{Viscous Burger's equation}
\begin{figure}[!h]
	\centerline{
		\includegraphics[scale=0.45]{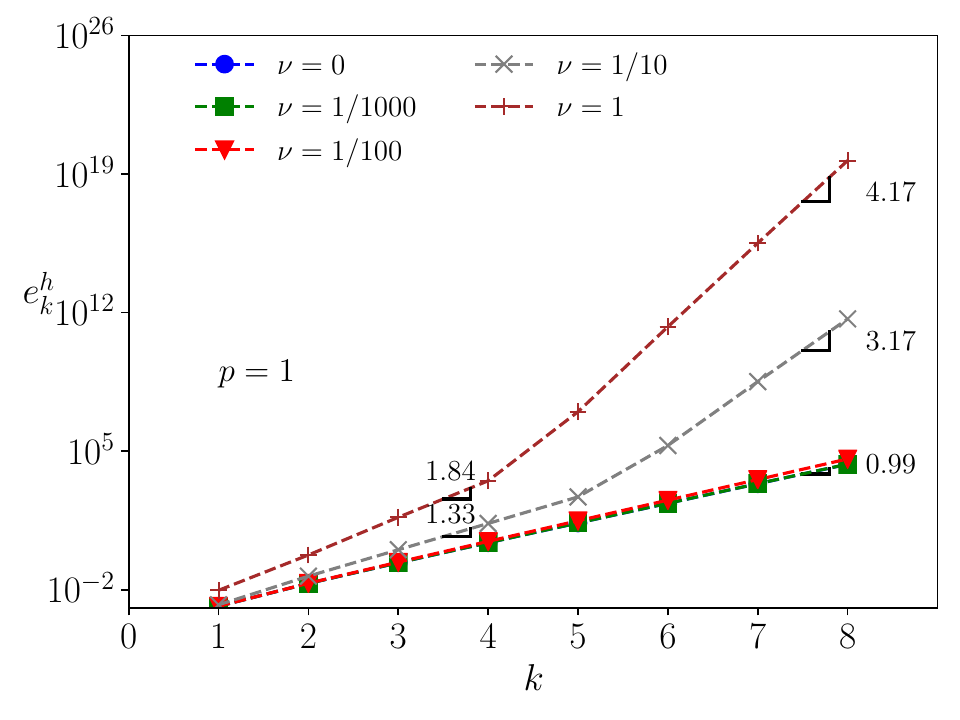}
		\includegraphics[scale=0.45]{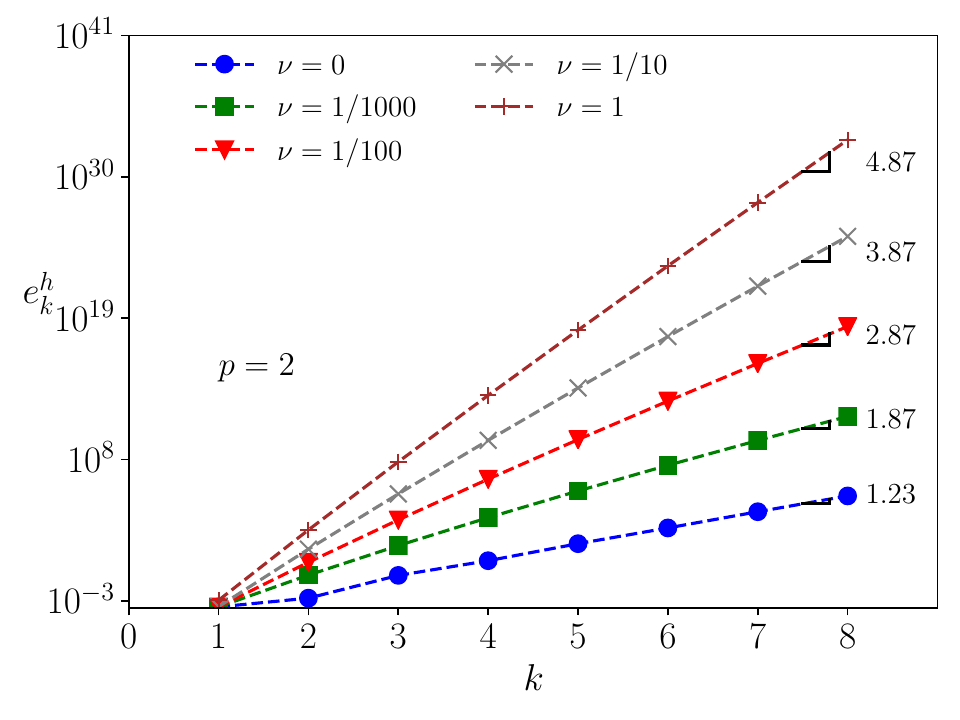}
	}
	\centerline{
		\includegraphics[scale=0.45]{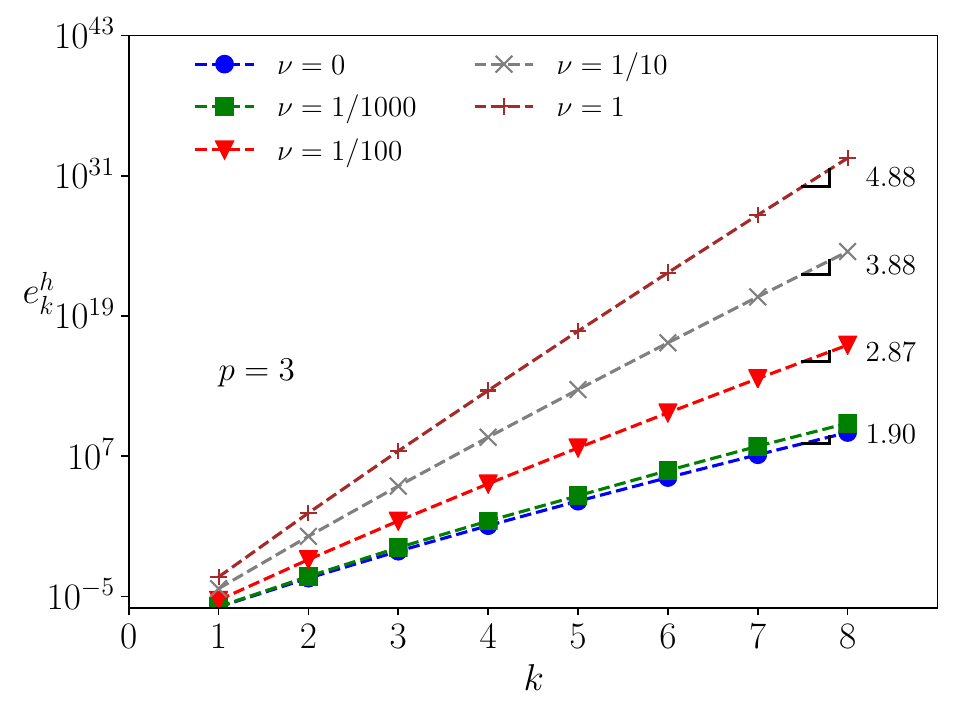}
		\includegraphics[scale=0.45]{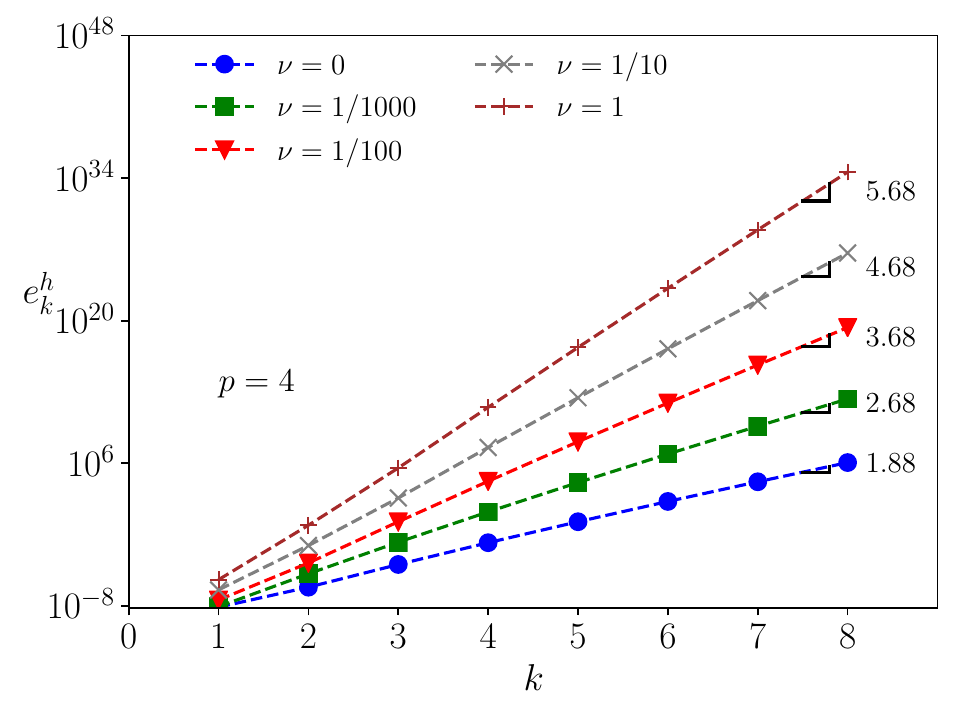}
	}
	\caption{Relative errors $\Me_k^{h,p}$ and their slopes as presented by formula \eqref{error_slope} for $h=1/100$, different degrees $p$ and different viscosity coefficients $\nu \in \{0,1/1000,1/100,1/10,1\}$}
	\label{fig11}
\end{figure}

We take now $\nu > 0$, we compute terms of the series without artificial diffusion for a fixed mesh size $h=1/100$ and different values of the viscosity $\nu$ and degree $p$. The error $\Me^h_k$ is evaluated and plotted in Figure \ref{fig11}.
We can see that in the case of $p=1$ and $0<\nu\leqslant1/100$, the error $\Me_k^h$ has an increasing rate of almost one in the logarithmic scale, which means that $\Me_k^h \simeq \Fc(h)\times10^k$. For greater values of viscosity, $\nu>1/100$, the error presents a steady slope for the first values of $k$: $\slope\simeq 1.33$ for $\nu=1/10$ and $k\leqslant5$, while $\slope \simeq 1.84$ for $\nu=1$ and $k\leqslant 4$, before jumping to a bigger slope that is equal to $\slope =3.17$ when $\nu=1/10$ and $\slope\simeq 4.17$ when $\nu=1$. We already have seen a such jump in the heat equation when $p=1$. For higher degrees $p\geqslant 2$, the errors also increase with $k$ and the viscosity $\nu$, having a fixed rate for a fixed viscosity. The slope of every line is showed in the plots. We can approximate the evolution of the error in function of $k$ as presented by Formula \eqref{error_slope}, thus the amplification error is just the power of $10$ by this slope $\slope(h,p)$.

Checking again the values of the slopes, we remark first that for a fixed \ac{FEMF}, the slopes increases by $1$ when $\nu$ is multiplied by ten. This shows that the error rate, \emph{i.e.} the slope $\slope$ depends globally on $\log_{10}(\nu \,\K^h)$, though it could be approximated by the logarithmic, in base 10, of the amplification factor given in Formula \eqref{ampl_fact_burg}, as $\HS{\K^h} \gg \HS{\D^h}$.

\begin{figure}[!h]
	\centerline{
		\includegraphics[scale=0.45]{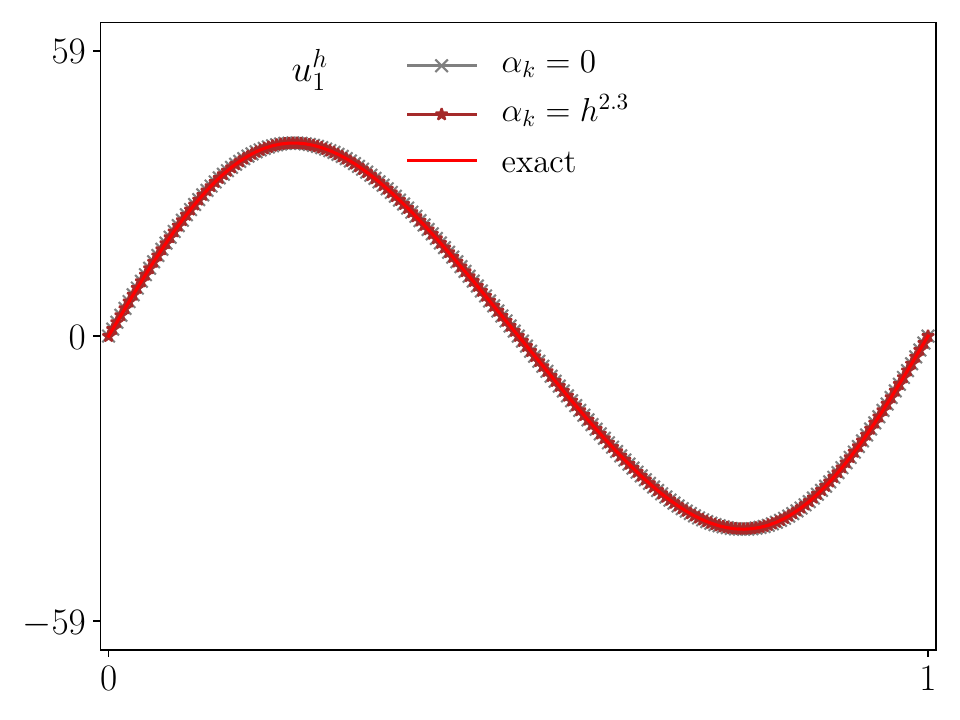}
		\includegraphics[scale=0.45]{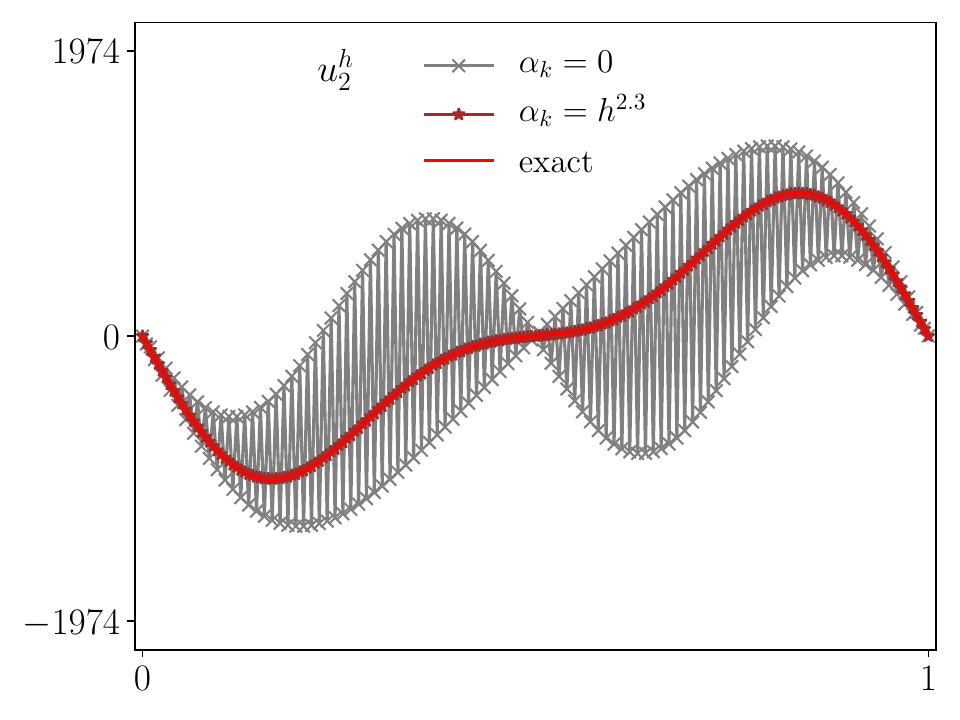}
	}
	\centerline{
		\includegraphics[scale=0.45]{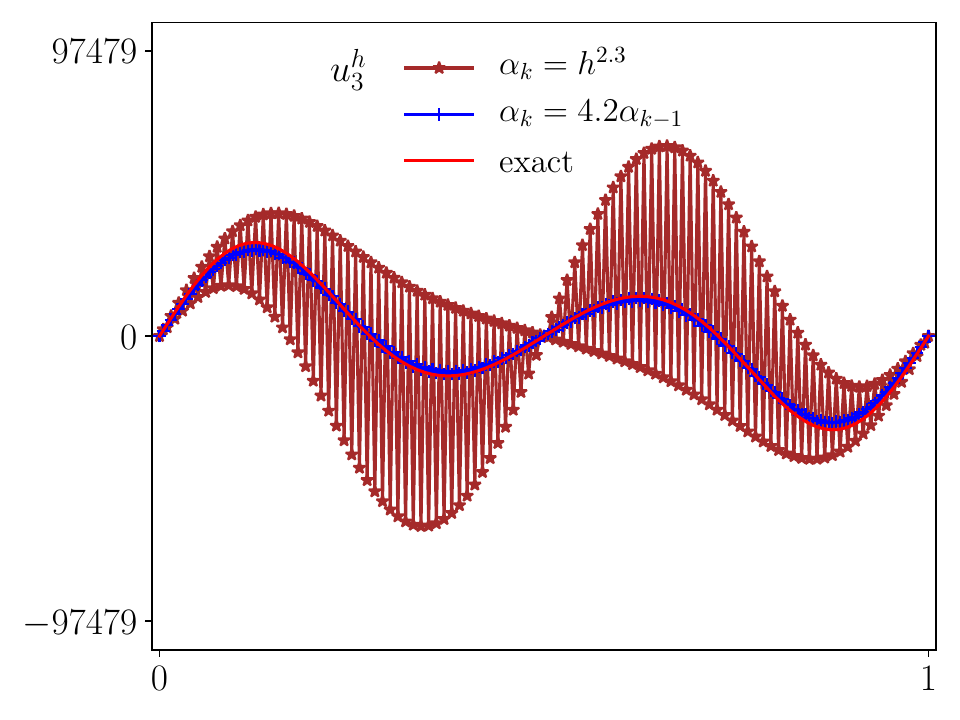}
		\includegraphics[scale=0.45]{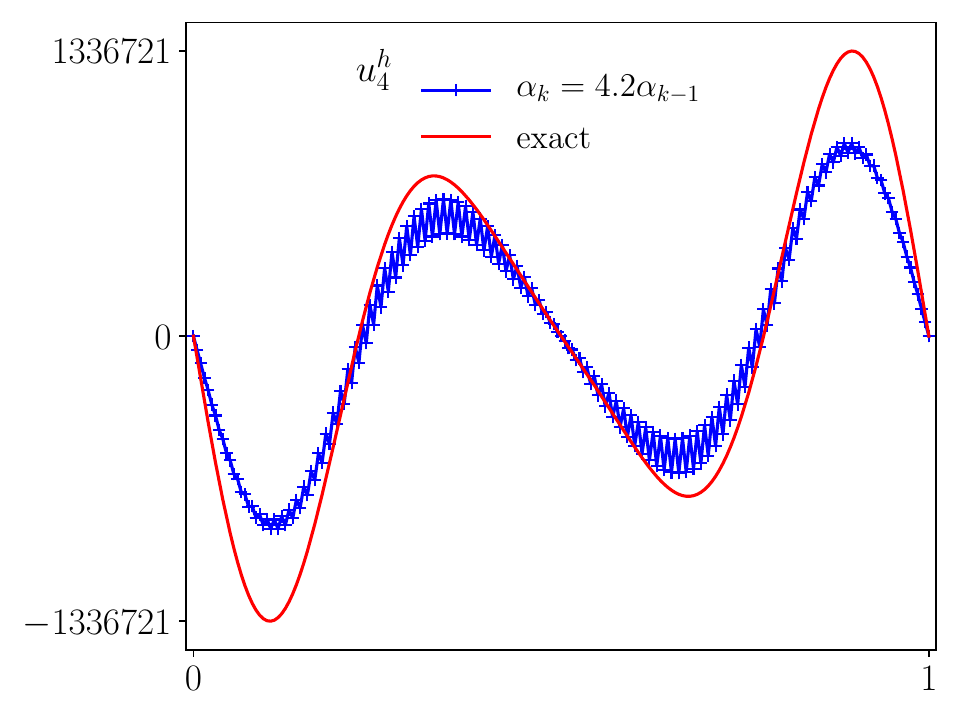}
	}
	\caption{Plots of the exact terms $\Mu_k, \,k\in\{ 1,2,3,4\}$ of the viscous Burger equation's series for $\nu=1$, their approximations in the space $\mathbb{V}_{h,p}$ for $h= 1/100$ and $p = 2$, with and without stabilization.}
	\label{fig12}
\end{figure}

We apply the stabilization technique presented in Section \ref{sec:alg}. Figure \ref{fig12} presents the plots of exact solutions of $u_k$ for $k\in\,\{1,2,3,4\}$ and their approximations in the \ac{FEMF} with $h=100$ and $p=2$. These approximations are obtained first without stabilization technique ($\alpha_k=0$), which approximates the first term with a high precision that is almost equal that one of the stabilization technique with fixed coefficient, $\alpha_k = h^{c}$, where $c\simeq 2.3$ is used from Table \ref{tab-heat-alpha-1d}. However, computing $u_2^h$ without the stabilization produces a second term that oscillates around the exact one, while the stabilization with fixed coefficient succeed in computing a reliable second term.
Nevertheless, it fails to produce the third term $u_3^h$ whilst the one with dependence stabilization coefficient, $\alpha_0 = h^{c}$ and $\alpha_k=R(h,p)\alpha_{k-1}$ with $R(h,p)\simeq 4.2$ is fetched from Table \ref{tab-heat-alpha-1d}, has produced a reliable third term that will be used in computing the fourth one. The last plot shows only the exact $u_4$ and the approximation with the stabilization technique having a dependence coefficient as it was the only one that succeed. This proves the importance of computing terms of the series using the proposed strategy.

To this end, we run the simulation: we compute the terms of the series with and without stabilization up to order $m=5$ and using several \ac{FEMF}. Then, the \ac{BPL} algorithm is used for approximating $\Mu^h(t)$ and advancing in time with different values of fixed time step $\Delta t^{\BPL}$. The integral residual errors in \eqref{Res} are computed and values are reported in Table \ref{tab8} for $T=0.5$. The mark x means that the relative simulation for the given time step and \ac{FEMF} ($h$ and $p$) does not reach the final time $T$ as the residual error \eqref{residual-error} explodes in times. We conclude from this table that applying the \ac{BPL} algorithm while computing terms of the series without the proposed stabilization inflicts on us to have a smaller time step to keep the stability criteria and reach the final instant. However, the proposed stabilization technique allows us to produce the same precision of the numerical approximation while having a larger time step; \emph{i.e.} for $p=2$ and $h=1/20$, the algorithm produces a residual error of $0.0194313$ with a time step of $\Delta t^{\BPL} = 5\times 10^{-2}$ and when employing the stabilization, however the time step should be smaller 200 times ($10^{-4}$) to produces a residual error of $0.017424$ without the stabilization.

\begin{sidewaystable}
	\caption{Residual errors \eqref{Res} for different features of simulation.}
	\label{tab8}
	\begin{center}
		\begin{tabular*}{\textheight}{@{\extracolsep\fill}cccc@{\extracolsep\fill}}
			\toprule
			& $\alpha_k=0$& $\alpha_k=h^c$ &  $\alpha_k=R(h,p)\alpha_{k-1}$    \\ \midrule
			$\Delta t^{\BPL} = 5\times10^{-2}$
			&
			\begin{tabular}{c|ccc}
				$h$ $\backslash$ $p$  & 1 & 2 & 3 \\ \midrule
				1/20 & ~~~~~x~~~~~ &~~~~~x~~~~~ &  ~~~~~x~~~~~\\
				1/50 & x & x &  x\\
				1/100  &x  &  x & x \\
				1/200 & x & x & x
			\end{tabular}
			&
			\begin{tabular}{ccc}
				1 & 2 & 3 \\ \midrule
				~~~~x~~~~ & ~~~~x~~~~ &  ~~~~~x~~~~~\\
				x & x &  x\\
				x  &  x & x \\
				x & x & x
			\end{tabular}
			&\begin{tabular}{ccc}
				1 & 2 & 3 \\ \midrule
				0.029671 & 0.0194313 &  0.017131\\
				0.018912 & x &  x\\
				x  &  x & x \\
				x & x & x
			\end{tabular}
			\\ \midrule
			$\Delta t^{\BPL} = 10^{-2}$
			&
			\begin{tabular}{c|ccc}
				$h$ $\backslash$ $p$  & 1 & 2 & 3 \\ \midrule
				1/20 & ~~~~~x~~~~~ &~~~~~x~~~~~ &  ~~~~~x~~~~~\\
				1/50 & x & x &  x\\
				1/100  &x  &  x & x \\
				1/200 & x & x & x
			\end{tabular}
			&
			\begin{tabular}{ccc}
				1 & 2 & 3 \\ \midrule
				0.02718 & 0.017171 &  ~~~~x~~~~\\
				0.01750 & x &  x\\
				x  &  x & x \\
				x & x & x
			\end{tabular}
			&\begin{tabular}{ccc}
				1 & 2 & 3 \\ \midrule
				0.02714 & 0.017159 &  0.01306\\
				0.01733 & 0.01094 &  x\\
				0.01225  &  x & x \\
				x & x & x
			\end{tabular}
			\\ \midrule
			$\Delta t^{\BPL} = 10^{-3}$
			&\begin{tabular}{c|ccc}
				$h$ $\backslash$ $p$ & 1 & 2 & 3 \\ \midrule
				1/20 &~ 0.02750~ & ~~~~~x~~~~~ & ~~~~x~~~~ \\
				1/50 & x & x &  x\\
				1/100 & x  &  x & x \\
				1/200 & x & x & x
			\end{tabular}
			&
			\begin{tabular}{ccc}
				1 & 2 & 3 \\ \midrule
				0.027244 & 0.017220& 0.013100 \\
				0.017284 & 0.010977 &  0.008315\\
				0.012294  &  0.007778 & x \\
				x & x  & x
			\end{tabular}
			&
			\begin{tabular}{ccc}
				1 & 2 & 3 \\ \midrule
				0.027246 & 0.017222  & 0.013101 \\
				0.017286 & 0.010977 & 0.008315 \\
				0.012295  &  0.007778 & 0.005884 \\
				0.008704&  0.005504 & x
			\end{tabular}
			\\ \midrule
			$\Delta t^{\BPL} = 10^{-4}$
			&\begin{tabular}{c|ccc}
				$h$ $\backslash$ $p$ & 1 & 2& 3 \\ \midrule
				1/20& 0.027516 & 0.017424 &0.013171   \\
				1/50 & 0.017420  & x &  x\\
				1/100 & 0.012320 &  x & x\\
				1/200 & x & x& x
			\end{tabular}
			&
			\begin{tabular}{ccc}
				1 & 2 & 3 \\ \midrule
				0.027263 & 0.017233  & 0.013108 \\
				0.017295 & 0.010983 & 0.008320 \\
				0.012301  &  0.007782 & 0.005887 \\
				0.008708 &  0.005507   & 0.004164
			\end{tabular}
			&
			\begin{tabular}{ccc}
				1 & 2 & 3 \\ \midrule
				0.027263 & 0.017233 &0.013108\\
				0.017295 & 0.010983 &0.008320 \\
				0.012301 & 0.007783 &0.005887\\
				0.008709&   0.005507 & 0.004164
			\end{tabular}
			\\
			\bottomrule
		\end{tabular*}
	\end{center}
\end{sidewaystable}

\section{Conclusions and perspectives}\label{sec-conc}

\subsection*{Conclusions.}
In this work, we present a new process to compute the terms of the series' solution type for a given \ac{PDE}. This process is based on a stabilization technique on the left hand side of the recurrence formula that generates these terms. The goal of this proposition is to provide stabilized terms in order to use the time series expansion in approximating the solution after applying one of the temporal integrator that is based on divergent series resummation techniques.

After presenting the ac{BPL} algorithm to integrate system of \ac{ODE}s and the continuation process to advance in time by building a sequences of approximations on set of points $t_n$, the problem in computing terms of the series in higher order \ac{FEMF} for Parabolic \ac{PDE}s has been presented. We have shown that terms of the series explode, though they are not usable if the precision of the \ac{FEMF} is refined by either increasing the mesh size or the degree of the \ac{FE}. We have shown also that using lumping techniques can help improve the stability of the errors but the results were not satisfying. An explanation of this explosion was linked to the condition number of the associated Mass matrix. The latter is a magnification factor of error propagation along the sequence of terms.

A modified system was proposed to reduce this magnification factor by adding on the left hand side of the recurrence formula an artificial diffusion weighted by a coefficient $\alpha_0$. This coefficient is sought as it minimizes the condition number of the new system and was considered depending on the \ac{FEMF} and the order $k$ of the computed term too. The link between the \acf{DMP} and the proposed stabilization technique was established, where the minization process contributes in satisfying the maximum principle int he discrete form. The stabilization coefficient $\alpha_k$ obeys a recurrence formula that is approximated in this by an empirical formula $\alpha_k = R(h,p)\times\alpha_{k-1}$. The ratio $R(h,p)$ is sought for every \ac{FEMF} so that it minimizes a new functional that depends on the amplification factor.

The proposed method was applied to the one dimensional heat and viscous Burger equations. Once the \ac{FEMF} is defined, the mass matrix and the stiff matrix are built and coefficient $\alpha_0$ and $R(h,p)$ are found. For a given initial condition for both equations, the terms are computed with and without stabilization, where the latter demonstrates how crucial is using it in providing stable and reliable terms of the series. Thought the employment of the proposed stabilization allows us to increase the time step by $200$ in the \ac{BPL} algorithm and having a residual error that is the same if the \ac{BPL} is applied without the stabilization.

\subsection*{Perspectives.}
Despite the improvement of computing stable and reliable terms of the series, we need to study the stability of the \ac{BPL} algorithm in producing numerical approximation for large time scale, \emph{i.e.} investigate the error propagation and the effect of the artificial diffusion in modifying the solution along the simulation process. Another question should be addressed about finding a priori-estimate of the approximation given by the \ac{BPL} integrator for a pre-defined used accuracy. This will help us (a) reducing the computational cost as we are evaluating the residual error for a sequences of time steps $\Delta t^{\BPL}$ until reaching a pre-defined user limit, and (b) improving the algorithm as the initial time step in the sequence is an estimation based on the partial sum. This estimation could be too small if the series is divergent.

In this work, we have used the pointwise \ac{Pa-ap} to prolongate analytically the series in the Borel-space. In the future, we need to check the effect of using Vector/Matrix Padé-type approximant in the case of one or two dimensional space, as we know how better is using vector-Padé type in approximating one-variable scalar functions.

\renewcommand{\nomname}{{NOMENCLATURE}}
\printnomenclature

\bmsection*{List of Abbreviations}
\begin{acronym}[TDMA]
\acro{PDE}{Partial Differential Equation}
\acro{DMP}{Discrete Maximum Principle}
\acro{TSE}{Time Series Expansion}
\acro{FE}{Finite Element}
\acro{FEM}{Finite Element Method}
\acro{ODE}{Ordinary Differential Equation}
\acro{BDF}{Backward Difference Formulas}
\acro{FDM}{Finite Difference Methods}
\acro{FVM}{Finite Volume Methods}
\acro{PGD}{Proper Generalized Decomposition}
\acro{FEMF}{FEM Framework}
\acro{NS}{Navier-Stokes}
\acro{BPL}{Borel-Padé-Laplace}
\acro{IVP}{Initial Value Problem}
\acro{LSM}{Level Set Method}
\acro{Pa-ap}{Padé approximant}
\acro{SVD}{Singular Value Decomposition}
\acro{GL}{Gau\ss-Laguerre}
\end{acronym}


\bmsection*{Acknowledgments}
This publication is based upon work supported by the Khalifa University of Science and Technology under Award No. FSU-2023-014.

\bmsection*{Financial disclosure}

None reported.

\bmsection*{Conflict of interest}

The authors declare no potential conflict of interests.

\bibliography{biblio}
\end{document}